\UseRawInputEncoding
\documentclass[11pt]{article}
\usepackage{amsmath}
\usepackage[numbers,sort&compress]{natbib}
\numberwithin{equation}{section}
\usepackage{amsfonts,amssymb,amsmath}
\usepackage{mathrsfs}
\usepackage{color}
\usepackage[colorlinks,bookmarksopen,bookmarksnumbered,citecolor=black, linkcolor=blue, urlcolor=black]{hyperref}
\newcommand{\beq}{\begin{equation}}
\newcommand{\enq}{\end{equation}}

\newtheorem{Theorem}{Theorem}[section]
\newtheorem{Lemma}[Theorem]{Lemma}
\newtheorem{Corollary}[Theorem]{Corollary}
\newtheorem{Definition}[Theorem]{Definition}
\newtheorem{Remark}[Theorem]{Remark}

\newcommand{\benu}{\begin{enumerate}}
\newcommand{\beqa}{\begin{eqnarray}}
\newcommand{\beqan}{\begin{eqnarray*}}
\newcommand{\eay}{\end{array}}
\newcommand{\edm}{\end{displaymath}}
\newcommand{\eenu}{\end{enumerate}}
\newcommand{\eeq}{\end{equation}}
\newcommand{\eeqa}{\end{eqnarray}}
\newcommand{\eeqan}{\end{eqnarray*}}

\newcommand{\br}{\begin{Remark}}
\newcommand{\er}{\end{Remark}}

\newcommand{\bqa}{\begin{eqnarray}}
\newcommand{\eqa}{\end{eqnarray}}
\newcommand{\bqw}{\begin{eqnarray*}}
\newcommand{\eqw}{\end{eqnarray*}}

\newcommand{\non}{\nonumber}
\newcommand{\bea}{\begin{array}{cc}}
\newcommand{\ena}{\end{array}}

\allowdisplaybreaks[4]
\begin{document}
\begin{center}

{\large \bf Existence and regularity of pullback attractors for nonclassical non-autonomous diffusion equations with delay}\\

\vspace{0.20in}Bin Yang$^{1}$ $\ $ Yuming Qin $^{2,\ast}$ $\ $ Alain Miranville $^{3}$ $\ $ Ke Wang $^{4}$\\
\end{center}
$^{1}$ College of Information Science and Technology, Donghua University, Shanghai, 201620, P. R. China.\\
$^{2}$ Department of  Mathematics, Institute for Nonlinear Science, Donghua University, Shanghai, 201620, P. R. China.\\
$^{3}$ Laboratoire de Math\'ematiques et Applications, Universit\'e de Poitiers, UMR CNRS 7348-SP2MI, Boulevard Marie et Pierre Curie-T\'el\'eport 2, F-86962, Chasseneuil Futuroscope Cedex, France.
 \vspace{3mm}

\begin{abstract}
In this paper, we consider the asymptotic behavior of weak solutions for non-autonomous diffusion equations with delay in time-dependent spaces when the nonlinear function $f$ is critical growth, the delay term $g(t, u_t)$ contains some hereditary
characteristics and the external force $h \in L_{l o c}^{2}\left(\mathbb{R} ; L^{2}(\Omega)\right)$. Firstly, we prove the well-posedness of solutions by using the Faedo-Galerkin approximation method. Then after a series of elaborate energy estimates and calculations, we establish the existence and regularity of pullback attractors in time-dependent spaces $C_{\mathcal{H}_{t}(\Omega)}$ and $C_{\mathcal{H}^{1}_{t}(\Omega)}$ respectively.
\end{abstract}

\hspace{4mm}{\bf Keywords:} Non-autonomous diffusion equations; Pullback attractors; Regularity.

\hspace{4mm}{\bf Mathematics Subject Classification (2020):} 35B40, 35B41, 35B65, 35K57.
\section{Introduction}
\setcounter{equation}{0}

\let\thefootnote\relax\footnote{*Corresponding author.}
\let\thefootnote\relax\footnote{E-mails: binyangdhu@163.com, Alain.Miranville@math.univ-poitiers.fr, yuming$\_$qin@hotmail.com.}
\quad
We shall investigate the following nonclassical non-autonomous diffusion equations with nonlocal term and delay
\begin{equation}
\left\{\begin{array}{ll}
\partial_{t}u-\varepsilon(t) \partial_{t}\Delta u-a(l(u)) \Delta u=f(u)+g(t, u_{t})+h(t) & \text { in } \Omega \times(\tau, \infty), \\
u(x,t)=0 & \text { on } \partial \Omega\times(\tau, \infty), \\
u(x, \tau+\theta)=\phi(x, \theta),  &\,\, x \in \Omega,\, \theta \in[-k, 0],
\end{array}\right.\label{1.1-3}
\end{equation}
in time-dependent spaces $C_{\mathcal{H}_{t}(\Omega)}$ and $C_{\mathcal{H}^{1}_{t}(\Omega)}$, where $k>0$ is the length of the delay effects, $\phi \in C\left([-k, 0] ; \mathcal{H}_{t}(\Omega)\right)$, $f$, $g$ and $h \in L_{l o c}^{2}\left(\mathbb{R} ; L^{2}(\Omega)\right)$ are nonlinear function, delay operator and external force function respectively. Besides, $\Omega  \subset \mathbb{R}^{n}\,(n \ge 3)$ is a bounded domain with smooth boundary $\partial\Omega$ and let $u_{t}(\theta)=u(t+\theta)$ and $u_t$ is defined in $[-k, 0]$ for any $t \in \mathbb R$.

Let the time-dependent function $\varepsilon(t) \in C^{1}(\mathbb{R})$ and satisfy
\begin{equation}
\lim _{t \rightarrow+\infty} \varepsilon(t)>\alpha>\frac{1}{2}
\label{1.2-3}
\end{equation}
and there exists a constant $L>0$ such that
\begin{equation}
\sup _{t \in \mathbb{R}}(|\varepsilon(t)|+|\varepsilon^{\prime}(t)|) \leq L.
\label{1.3-3}
\end{equation}

Besides, $a(l(u)) \in C\left(\mathbb{R}; \mathbb{R}^{+}\right)$ is a nonlocal diffusion function, where $l(u): L^{2}(\Omega) \to \mathbb R$ is a continuous linear functional acting on $u$ and $l(u)=l_i(u)=\int_{\Omega} i(x) u(x) d x$ for some $i \in L^{2}(\Omega)$. When $\varepsilon(t)$ is decreasing, $a(\cdot)$ satisfies
\begin{equation}\label{1.4-3}
m \leqslant a(s) \leqslant M
\end{equation}
and when $\varepsilon(t)$ is increasing, $a(\cdot)$ fulfills
\begin{equation}\label{1.5-3}
m+L \leqslant a(s) \leqslant M.
\end{equation}

Furthermore, suppose $f$ is Lipschitz continuous and satisfies $f(0)=0$,
\begin{equation}\label{1.6-3}
\limsup _{|u| \rightarrow \infty} \frac{f(u)}{u}<\lambda_1
\end{equation}
and
\begin{equation}\label{1.7-3}
\left|f^{\prime}(u)\right| \leq C\left(1+|u|^{p}\right),
\end{equation}
where $p=\frac{4}{n-2}$ and $\lambda_1>0$ is the first eigenvalue of $-\triangle$. Moreover, $f$ admits the decomposition $f=f_0+f_1$ with $f_0, f_1 \in C(\mathbb{R}, \mathbb{R})$,
where $f_0$ satisfies
\begin{equation}\label{1.8-3}
\left|f_0(u)\right| \leqslant C\left(|u|+|u|^{p+1}\right)
\end{equation}
and
\begin{equation}\label{1.9-3}
f_0(u) u \leq -l
\end{equation}
for any constant $l>0$ and $f_1$ satisfies $(\ref{1.6-3})$
and
\begin{equation}\label{1.10-3}
\left|f_1(u)\right| \leqslant C\left(1+|u|^\gamma\right)
\end{equation}
for any $0<\gamma <\frac{n+2}{n-2}$.

In addition, the delay $g: \mathbb{R} \times C_{L^{2}(\Omega)} \rightarrow L^{2}(\Omega)$ contains some hereditary
characteristics and satisfies
\begin{equation}
the\,\,function\,\,\mathbb{R}\,\,\ni\,\,t\,\,\mapsto\,\,g(t, u) \in L^{2}(\Omega)\,\,is\,\, measurable,\,\, \forall \, u\,\in\,C_{L^{2}(\Omega)},
\label{1.11-3}
\end{equation}
\begin{equation}
g(t, 0)=0, \,\, \forall\,t\,\in\,\mathbb{R},
\label{1.12-3}
\end{equation}
and there exists a constant $C_{g}>0$ such that
\begin{equation}
\|g(t, u_{1})-g(t, u_{2})\|^{2} \leq C_{g}\|u_{1}-u_{2}\|^{2}_{C_{L^{2}(\Omega)}},\,\, \forall \, u_{1}, u_{2} \in C_{L^{2}(\Omega)}.
\label{1.13-3}
\end{equation}

For the sake of simplicity, in this paper the norm and inner product of $L^{2}(\Omega)$ are written as $\|\cdot\|$ and $(\cdot, \cdot)$ respectively. Our phase space is the time-dependent space $C_{\mathcal{H}_{t}(\Omega)}$ and it is equipped with the norm
\begin{equation}\label{1.14-3}
\|u\|_{C_{\mathcal{H}_t(\Omega)}}^2=\|u\|_{C_{L^2(\Omega)}}^2+\left|\varepsilon_t\right|\|\nabla u\|_{C_{L^2(\Omega)}}^2,
\end{equation}
where $\|u\|_{C_{L^{2}(\Omega)}}=\max\limits _{\theta \in[-k, 0]}\|u(t+\theta)\|$ for any $t\in\mathbb R$ and $\left|\varepsilon_{t}\right|$ is the absolute value of $\varepsilon(t+\theta)$. Moreover, the time-dependent space $C_{\mathcal{H}_t^1(\Omega)}$ is more regular than $C_{\mathcal{H}_{t}(\Omega)}$ and it is endowed with the norm
\begin{equation}\label{1.15-3}
\|u\|_{C_{\mathcal{H}_t^1(\Omega)}}^2=\|\nabla u\|_{C_{L^2(\Omega)}}^2+\left|\varepsilon_t\right|\|\Delta u\|_{C_{L^2(\Omega)}}^2 .
\end{equation}
Similarly, define the norm of the time-dependent space $C_{\mathcal{H}_t^1(\Omega),\sigma}$ as
\begin{equation}\label{1.16-3}
\|u\|_{C_{\mathcal{H}^1_t(\Omega), \sigma}}^2=\|A^{\frac{\sigma}{2}} u\|_{C_{L^2(\Omega)}}^2+\left|\varepsilon_t\right|\|A^{\frac{1+\sigma}{2}} u\|_{C_{L^2(\Omega)}}^2,
\end{equation}
where the operator $A=-\Delta$ and $0<\sigma<\min \left\{\frac{1}{3}, \frac{n+2-(n-2) \gamma}{2}\right\}$ and the range of $\gamma$ is the same as in $(\ref{1.10-3})$.

In 2011, Plinio, Duane and Temam \cite{pdt.3} first proposed the time-dependent space $\mathcal H_{t}(\Omega)$ and defined its norm as $\|\cdot\|^{2}+\varepsilon(t)\|\nabla \cdot\|^{2}$, where $\varepsilon(t)$ must be decreasing and satisfies $\lim\limits _{t \rightarrow+\infty} \varepsilon(t)=0$. They called $\mathscr{A}=\left\{\mathscr{A}(t) \subset \mathcal H_t\right\}_{t \in \mathbb{R}}$ is a time-dependent global attractor, if it fulfills
(i) $\mathscr{A}(t)$ is compact; (ii) $S(t, s) \mathscr{A}(s)=\mathscr{A}(t)$, for every $s \leq t$;
(iii) $\lim \limits_{s \rightarrow-\infty} \operatorname{dist}_{{\mathcal H}_t}(S(t, s) \mathcal{B}(s), \mathscr{A}(t))=0$, for every pullback-bounded family $\mathcal{B}$ and every $t \in \mathbb{R}$.
Moreover, if property (ii) holds uniformly with respect to $t \in \mathbb{R}$, $\mathscr{A}$ is a uniform time-dependent global attractor.
It is easy to derive if $\mathcal{B}$ is a family of all pullback-bounded family in $\mathcal H_{t}(\Omega)$, then $\mathscr A$ is a pullback $\mathcal D$-attractor.
Then Conti, Pata and Temam \cite{cpt.3} redefined this definition by no longer requiring invariance, and proved $\mathcal A$ should admits (i) and ${\mathscr{A}}$ is pullback attracting, i.e., it is uniformly bounded and satisfies $\lim\limits_{\tau \rightarrow-\infty} dist \left(U(t, \tau) B_\tau, \mathscr{A}_t\right)=0$, where $\widehat{B}=\left\{B_t\right\}_{t \in \mathbb{R}}$ is a uniformly bounded family and $t \in \mathbb{R}$ is fixed. Besides, they also proved $\mathcal A$ exists and it is unique if and only if the process is asymptotically compact.

Inspired by above works, many scholars have focused on the time-dependent global attractors of various systems in $\mathcal H_{t}(\Omega)$. Conti and Pata \cite{CP} obtained the existence of the time-dependent global attractors for wave equation $\varepsilon u_{tt}+\alpha  u_t-\Delta{u}+f(u)=g(x)$. Meng, Yang and Zhong \cite{myz} studied wave equation $\varepsilon(t) u_{t t}+g\left(u_{t}\right)-\Delta u+\varphi(u)=f(x)$ and established a necessary and sufficient condition for the existence of the time-dependent global attractors. Ding and Liu \cite{DL} considered diffusion equation $u_t - \varepsilon(t)\Delta{u_t}-\Delta{u}+f(u)=g(x)$. Besides, Ma, Wang and Xu \cite{mwx} learned the existence and regularity of the time-dependent global attractors for reaction-diffusion equation $u_{t}-\varepsilon(t) \triangle u_{t}-\Delta u+\lambda u=f(u)+g(x)$. Zhu, Xie and Zhou \cite{zxz} verified the existence of the time-dependent global attractors for reaction-diffusion equation $u_{t}-\varepsilon(t) \Delta u_{t}-\Delta u+f(u)=g(x)$. Moreover, Meng, Wu and Zhao \cite{mwz} investigated the time-dependent global attractors of extensible Berger equation $\varepsilon(t) u_{t t}+\Delta^{2} u-\left(Q+\int_{\Omega}|\nabla u|^{2} d x\right) \Delta u+g\left(u_{t}\right)+\varphi(u)=f(x)$, where $Q$ is a plane internal force applied on the plate. Wang, Hu and Gao \cite{whg.3} considered undamped abstract evolution equation $\varepsilon(t) u_{t t}+k(0) A^\theta u+\int_0^{\infty} k^{\prime}(s) A^\theta u(t-s) \mathrm{d} s+f(u)=g(x)$, where $\theta \in\left({2 n/(n+2)}, n/2\right)$ and $k(\cdot)$ is a memory kernel. In addition, there are some other relevant works, see \cite{hrz.3, ML, MWL}, etc. On the other hand, please refer to our paper \cite{qy3.3}, where the relevant papers to problem (\ref{1.1-3}) are comprehensively introduced.

Below we will introduce difficulties and our innovations. Motivated by Conti, Pata and Temam \cite{cpt.3}, Sun and Yang \cite{sy.3} and Zhu and Sun \cite{zs2}, we have obtained our main results. Compared our phase spaces $C_{\mathcal{H}_{t}(\Omega)}$ and $C_{\mathcal{H}^{1}_{t}(\Omega)}$ with $\mathcal H_t(\Omega)$ in \cite{cpt.3}, it can be found that we redefine the limit of the time-dependent term $\varepsilon(t)$ is equal to $\alpha$ when $t \rightarrow+\infty$, and we also take into account the longtime behavior of weak solutions to problem (\ref{1.1-3}) when $\varepsilon(t)$ is increasing. Meanwhile, due to the existence of nonlocal term $a(l(u))$ and delay $g(t,u_t)$, problem (\ref{1.1-3}) can be widely used in daily life.

Furthermore, the estimates of the existence and uniqueness of weak solutions are omitted in many papers. However, we make detailed additions to the proofs of this part, which will contribute to better understanding why we can only deduce the existence of weak solutions to problem (\ref{1.1-3}) but not strong solutions. In addition, some estimates here can also be used in the proofs of the existence of absorbing sets, asymptotically compactness of process and regularity of pullback attractors.

Moreover, our assumptions for the nonlinear function $f$ are weaker than those in \cite{zs2}. It is worth mentioning that the hypothesis of parameter $\sigma$ is $0<\sigma<\min \left\{\frac{1}{4}, \frac{n+2-(n-2) \gamma}{2}\right\}$ with $0<\gamma<\frac{n+2}{n-2}$ and $n \ge3$ in \cite{sy.3}, and it is directly quoted in some papers, such as \cite{zs2}. After we calculated their works in detail, it turns out that the range of $\sigma$ is set by some embedding inequalities, which are crucial to getting their results. Besides, we find that the  largest range of $\sigma$ is $0<\sigma<\min \left\{\frac{1}{3}, \frac{n+2-(n-2) \gamma}{2}\right\}$ and we will prove it in Lemma $\ref{lem4.4-3}$. As for papers that directly refer to the assumptions of $\sigma$ in \cite{sy.3}, their results are still valid under our hypothesis.

This paper is organized as follows. In $\S 2$, we shall recall some definitions and lemmas. Then, in $\S 3$ we verify the existence and uniqueness of solutions to problem $(\ref{1.1-3})$ in $C_{\mathcal{H}_{t}(\Omega)}$. Later on, in $\S 4$ we obtain the existence of pullback attractor $\mathcal A$ for the process ${\{ U(t,\tau )\} _{t \ge \tau }}$ in $C_{\mathcal{H}_{t}(\Omega)}$ and preliminarily verify the regularity of solutions in $C_{\mathcal{H}_t^1(\Omega),\sigma}$. Finally, in $\S 5$ we prove the regularity of solutions to problem $(\ref{1.1-3})$ in $C_{\mathcal{H}_t^1(\Omega)}$.

\section{\large Preliminaries}
In this section, in order to learn the asymptotic behavior of solutions in time-dependent spaces, we shall review some basic knowledge of pullback attractors and spaces.

First of all, suppose $\left\{X_t\right\}_{t \in \mathbb{R}}$ is a family of normed spaces.
\begin{Definition} {\rm(\cite{r})}\label{def2.1-3}
The Hausdorff semi-distance between two nonempty  sets $X_{a}, X_{b} \subset X_{t}$ is defined by
$$
dist_{X_{t}}(X_{a},X_{b})=\sup _{x_1 \in X_{a}} \inf _{x_2 \in X_{b}}\|x_1-x_2\|_{X_{t}} \, .
$$
\end{Definition}

\begin{Definition} {\rm(\cite{r})}\label{def2.2-3}
The set $\bar{\mathscr B}_{X_t}(0, R)$ is called a closed sphere centered at the origin and of radius $R$, if it satisfies
$$
\bar{\mathscr B}_{X_t}(0, R)=\left\{u \in X_t:\|u\|_{X_t} \leq R\right\}.
$$
\end{Definition}

\begin{Definition} {\rm(\cite{bcl.3, clr})} \label{def2.3-3}
A process in $\left\{X_t\right\}_{t \in \mathbb{R}}$ is the family $\{U(t, \tau)\}_{t \geqslant \tau}$ of mapping $U(t, \tau): X_{\tau} \rightarrow X_{t}$ that satisfies

(i) $U(\tau, \tau)=Id$ is the identity operator in $X_{\tau}$;

(ii) $U(t, s) U(s, \tau)=U(t, \tau)$ for any $t \geq s \ge \tau.$

\end{Definition}

\begin{Definition} {\rm(\cite{chm})}\label{def2.4-3}
For any constant $\tilde\epsilon>0$, if $\mathcal D$ is a nonempty class of all families of $\widehat{D}=\left\{D(t)\right\}_{t \in \mathbb{R}} \subset \Gamma(X_{t})$ such that
$$
\lim _{\tau \rightarrow-\infty}\left(e^{\tilde{\epsilon} \tau} \sup _{u \in D(\tau)}\|u\|_{X_{t}}^{2}\right)=0,
$$
where $\Gamma(X_{t})$ denotes a family of all nonempty subsets of $\left\{X_{t}\right\}_{t \in \mathbb{R}},$ then $\mathcal D$ is a tempered universe in $\Gamma(X_{t})$.
\end{Definition}

\begin{Definition} {\rm(\cite{psz,zs2})}\label{def2.5-3}
A family ${\widehat D_{0}=\left\{D_0(t)\right\}_{t \in \mathbb{R}}} \subset \Gamma(X_{t})$ is called pullback $\mathcal D$-absorbing for the process ${\{ U(t,\tau )\} _{t \ge \tau }}$, if for any $t \in \mathbb{R}$ and ${\widehat D} \in {{\cal D}}$, there exists a ${\tau _0} = {\tau _0}(t,{\widehat D}) < t$ such that
$
U(t, \tau) D(\tau) \subset D_{0}(t)
$
for any $\tau \leq \tau_{0}(t, \widehat{D}) \leq t \in \mathbb R$.
\end{Definition}

\begin{Definition} {\rm(\cite{Evans})}\label{def2.6-3}
Let $B \subset\left\{X_t\right\}_{t \in \mathbb R}$ be a bounded set, the noncompactness measure $k$ of $B$ is defined as
$$
k(B)=\inf \left\{\delta>0 \mid B \text { can be covered by a finite number of\,\,} d \text {-neighborhoods with} \, d \leqslant \delta \right\}.
$$
\end{Definition}

\begin{Lemma} {\rm(\cite{Evans})}\label{lem2.7-3}
Assume the sets $A_0, A_1$ and $A_2$ are bounded in $\left\{X_t\right\}_{t \in \mathbb{R}}$, then their noncompactness measures satisfies the following properties

i) $k\left(A_0\right)=0 \Leftrightarrow k\left(\mathcal N\left(A_0, \epsilon_0\right)\right) \leqslant 2 \epsilon_0$ $\Leftrightarrow \bar{A}_0$ is compact,

ii) $k\left(A_1+A_2\right) \leq k\left(A_1\right)+k\left(A_2\right)$,

iii) if $A_1 \subseteq A_2$, then $k\left(A_1\right) \leqslant k\left(A_2\right)$,

iv) $k\left(A_1, A_2\right) \leqslant \max \left\{k\left(A_1\right), k\left(A_2\right)\right\}$,

v) $k(\bar{A})=k(A)$,

\noindent where the symbol $\Leftrightarrow$ stands for if and only if, and $\mathcal N(A_0, \epsilon_0)$ denotes the neighborhood of $A_0$ whose radius is $\epsilon_0$.
\end{Lemma}

\begin{Definition} {\rm(\cite{bcl.3, clr})}\label{def2.8-3}
The process $\{U(t, \tau)\}_{t \geqslant \tau}$ is called pullback $\mathcal D$-$\omega$-limit compact, if for any constant $\tilde{\epsilon}>0$ and any bounded $\widehat{\mathcal D} \subset \mathcal D$, there exists a $\tilde{\tau}<t$ depending on $\tilde{\epsilon}$ and $\widehat{D}$ such that
$$
k\left(\bigcup_{\tau \leqslant \tilde{\tau}} U(t, \tau) D(\tau)\right) \leqslant \tilde{\epsilon}.
$$
\end{Definition}

\begin{Lemma} {\rm(\cite{bcl.3, clr})}\label{lem2.9-3}
Assume the process $\{U(t, \tau)\}_{t \geqslant \tau}$ can be decomposed as $U(t, \tau)=U_a(t, \tau)+U_b(t, \tau)$ for any fix $t \in \mathbb{R}$ and $\widehat{D}=\{D(t)\}_{t \in \mathbb{R}}$ is a pullback absorbing set for $\{U(t, \tau)\}_{t \geqslant \tau}$. Furthermore, if every subsequence of $U_a(t, \tau) D(\tau)$ is a Cauchy sequence for any fixed $\tau \leqslant t$ and $\lim\limits _{\tau \rightarrow-\infty}\left\|U_b(t, \tau) D(\tau)\right\|_{X_t}=0$, then $\{U(t, \tau)\}_{t \geqslant \tau}$ is pullback $\mathcal D$-$\omega$-limit compact in $\left\{X_t\right\}_{t \in \mathbb R}$.
\end{Lemma}

Furthermore, we will introduce a new space and its properties, which will be conducive to prove the asymptotic compactness and the regularity.

\begin{Definition} {\rm(\cite{ps.2})}\label{def2.10-3}
The space $D(A^{\frac{s}{2}})$ with $A=-\Delta$ and $s \in \mathbb{R}$ is a Hilbert space and its inner product and norm are defined as $(A^{\frac{s}{2}}\cdot, A^{\frac{s}{2}} \cdot)$ and $\|A^{\frac{s}{2}} \cdot\|$, respectively.
\end{Definition}

\begin{Lemma} {\rm(\cite{ps.2})}\label{lem2.11-3}
The properties of $D(A^{\frac{s}{2}})$ are as follows:

i) The embedding $D(A^{\frac{s_1}{2}}) \hookrightarrow D(A^{\frac{s_2}{2}})$ is compact for any $s_1>s_2$.

ii) The embedding $D(A^{\frac{s}{2}}) \hookrightarrow L^{\frac{2 n}{n-2 s}}\left(\Omega\right)$ is continuous for any $s \in\left[0, \frac{n}{2}\right)$.

iii) If $s_0>s_1>s_2$, then there exists constants $\epsilon_{1}>0$ and $C(\epsilon_{1})>0$ depending on $s_0$, $s_1$ and $s_2$ such that
$$
\|A^{\frac{s_1}{2}} \cdot\| \leq \epsilon_{1}\|A^{\frac{s_0}{2}} \cdot\|+C(\epsilon_{1})\|A^{\frac{s_2}{2}} \cdot\|.
$$

iv) Assume that $s_1, s_2 \in(0,1)$ and $u \in D(A^{\frac{s_1}{2}}) \cap D(A^{\frac{s_2}{2}})$, then there exists constants $\theta \in(0,1)$ and $C(\theta)>0$ such that
$$
\|\cdot\|_{D\left(A^{ \frac{(1-\theta) s_1+\theta s_2}{2}}\right)} \le C(\theta)\|\cdot\|_{D\left(A^{\frac{s_1}{2}}\right)}^{1-\theta}\|\cdot\|_{D\left(A^{\frac{s_2}{2}}\right)}^\theta.
$$
\end{Lemma}

\begin{Definition} {\rm(\cite{cv})}\label{def2.12-3}
The space $L_{b}^2\left(\mathbb{R};L^2(\Omega)\right)$ consists of all translation bounded functions in the space $L_{loc}^2\left(\mathbb{R};L^2(\Omega)\right)$, which is defined as
$$
L_b^2\left(\mathbb{R} ; L^2(\Omega)\right):=\left\{k \in L_{l o c}^2\left(\mathbb{R} ; L^2(\Omega)\right): \sup _{t \in \mathbb{R}} \int_t^{t+1} \int_{\Omega}|k(x, s)|^2 d x d s<+\infty\right\}.
$$
\end{Definition}

\begin{Definition} {\rm(\cite{bcl.3, clr})}\label{def2.13-3}
$A$ process $\{U(t, \tau)\}_{t \geq \tau}$ is called pullback $\mathcal D$-asymptotically compact in $\{X_t\}_{t \in \mathbb{R}}$, if for any $\widehat{D} \subset D$, any sequence $\left\{\tau_n\right\}_{n \in \mathbb N^{+}}\subset(-\infty, t]$ and any sequence $\left\{x_n\right\}_{n \in \mathbb N^{+}} \subset D\left(\tau_n\right) \subset X_t$, the sequence $\left\{U(t, \tau) x_n\right\}_{n \in \mathbb N^{+}}$is relatively compact in $\{X_t\}_{t \in \mathbb{R}}$.
\end{Definition}
\begin{Definition} {\rm(\cite{bcl.3, clr})}\label{def2.14-3}
A family ${\mathcal{A}}=\{\mathcal{A}(t)\}_{t \in \mathbb{R}} \subset \Gamma\left(X_{t}\right)$
is called a pullback $\mathcal D$-attractor for the process ${\{ U(t,\tau )\} _{t \ge \tau }}$ in $\left\{X_{t}\right\}_{t \in \mathbb{R}}$, if it satisfies the following properties

(i) $\mathcal A(t)$ is compact in $X_{t}$.

(ii) ${\mathcal A(t)}$ is pullback $\mathcal D$-attracting in $X_{t}$, i.e.,
$$
\lim _{\tau \rightarrow-\infty} {dist}_{X_{t}}\left(U(t, \tau) D(\tau), \mathcal A(t)\right)=0,
$$
for any ${\widehat D} \in {{\cal D}}$.

(iii) ${\mathcal{A}}$ is invariant, i.e., $U(t, \tau)\mathcal{A}(\tau)={\mathcal{A}(t)}$ for any $\tau \leq t \in \mathbb R$.
\end{Definition}

\begin{Lemma}{\rm(\cite{bcl.3, clr})}\label{lem2.15-3}
Let $B \subset X$ is bounded, then the process ${\{ U(t,\tau )\} _{t \ge \tau }}$ has a unique pullback $\mathcal D$-attractor $\mathcal A=\{\mathcal A(t)\}_{t \in \mathbb{R}}$ with
$$
\mathcal{A}(t)=\omega(B, t)=\bigcap_{\tau_1 \leq t} \overline{\bigcup_{\tau \leq \tau_1} U(t, \tau) B(\tau)},
$$
if and only if ${\{ U(t,\tau )\} _{t \ge \tau }}$ has a pullback $\mathcal D$-absorbing set and $\{U(t, \tau)\}_{t \geq \tau}$ is pullback $\mathcal D$-$\omega$-limit compact in $\left\{X_t\right\}_{t \in \mathbb{R}}$.
\end{Lemma}

\begin{Lemma} {\rm(\cite{Q1, Q2, Q3})}\label{lem2.16-3}
Let $u(t)$, $v(t)$, and $w(t)$ be real functions defined on interval $[a, b]$, where $u(t)$ is non-negative and Lebesgue integrable, $v(t)$ is absolutely continuous, $w(t)$ is continuous, which satisfy
$$
w(t) \leq v(t)+\int_a^t u(s) w(s) d s,\,\,\forall\,a \leq t \leq b,
$$
then the following inequality holds
$$
w(t) \leq v(a) e^{\int_a^t u(s) d s}+\int_a^t e^{\int_s^t u(\tau) d \tau} \cdot \frac{d v}{d s} d s
$$
for any $t \in [a, b]$.
\end{Lemma}

\section{Existence and uniqueness of weak solutions}
In this section, we will first define weak solutions to problem $(\ref{1.1-3})$ and then prove its existence by the energy method.

\begin{Definition}\label{def3.1-3}
Under the assumptions of $\varepsilon(t)$, $a(l(u))$, $f$, $g$ and $h$ in $\S 1$, $u$ is called a weak solution to problem $(\ref{1.1-3})$, if $u \in C([\tau-k, T]; \mathcal{H}_{t}(\Omega)) \cap L^{2}(\tau, T ; H_0^1(\Omega))$ with $u(t)=\phi(t-\tau)$ for any $\tau < T \in \mathbb R$ and $t \in[\tau-k, \tau]$, and satisfies
\begin{equation}
\frac{d}{d t}((u(t), \zeta)+\varepsilon(t)(\nabla u(t), \nabla \zeta))+(2 a(l(u))-\varepsilon^{\prime}(t))(\nabla u(t), \nabla \zeta) \\
= 2(f(u(t))+g(t, u_{t})+ h(t), \zeta),
\label{3.1-3}
\end{equation}
for any test function $\zeta \in H_{0}^{1}(\Omega)$.
\end{Definition}

\begin{Corollary}
If $u$ is a weak solution to problem $(\ref{1.1-3})$, then it satisfies
\begin{equation}
\begin{array}{l}
\|u(t)\|^{2}+\varepsilon(t)\|\nabla u(t)\|^{2}+\int_{s}^{t}\left(2 a(l(u))-\varepsilon^{\prime}(r)\right)\|\nabla u(r)\|^{2} d r \\
=\|u(s)\|^{2}+\varepsilon(s)\|\nabla u(s)\|^{2}+2 \int_{s}^{t}(f(u(r))+g(t, u_{r})+h(r), u(r)) dr,
\end{array}
\label{3.2-3}
\end{equation}
for any $s\in[\tau, t]$.
\end{Corollary}

Next, by using the Faedo-Galerkin approximation method, we can prove the following theorem about the existence and uniqueness of weak solutions to problem $(\ref{1.1-3})$.

\begin{Theorem}\label{th3.3-3}
Under the assumptions of $\varepsilon(t)$, $a(l(u))$, $f$, $g$ and $h$ in $\S 1$, assume $\phi \in C_{\mathcal H_{t}(\Omega)}$ is given, then there exists a weak solution $u(\cdot)=u(\cdot, \tau ; \phi)$ to problem $(\ref{1.1-3})$  in the time-dependent space $C_{\mathcal{H}_{t}(\Omega)}$, which satisfies $u \in C\left([\tau-k, T] ; \mathcal{H}_t(\Omega)\right)$ and $\partial_t u \in L^2\left(\tau, T ; \mathcal{H}_t(\Omega)\right)$
for any $\tau <T \in \mathbb R$. Moreover, the weak solution $u$ depends continuously on its initial value.
\end{Theorem}
$\mathbf{Proof.}$ Let $\left\{ {{e_j}} \right\}_{j \in \mathbb N^+} $ be a basis of $H^{2}(\Omega) \cap H_{0}^{1}(\Omega)$ and orthonormal in $L^{2}(\Omega)$. For the sake of using the Faedo-Galerkin approximation method, we need to find an approximate sequence $u_{i}(t, x)=\sum\limits_{j=1}^{i} \chi_{i, j}(t) \omega_{j}(x)$ with $i, j \in \mathbb N^+$ that satisfies
\begin{equation}
\left\{ {\begin{array}{*{20}{l}}
\frac{d}{dt}(({u_i}(t),{e _j}) + \varepsilon (t)(\nabla {u_i}(t),\nabla {e _j})) + ( {2a(l({u_i})) - {\varepsilon ^\prime }(t)} )(\nabla {u_i}(t),\nabla {e _j})\\
{ = 2(f({u_i}(t)),{e _j}) +2(g(t, u_{t,i}),{e _j})+ 2\left({h(t),{e _j}} \right), \,\, \forall\,\, t \in [\tau , + \infty ),}\,1 \leq j \leq i,\\
{u_{\tau,i}=\phi,}
\end{array}} \right.
\label{3.3-3}
\end{equation}
where $\phi \in C([-k, 0]; \operatorname{span}\{e_{j}\}_{j=1}^{n}), u_{t, i}=u_{i}(t+\theta)$ and $u_{\tau,i}=u_{i}(\tau+\theta)$ with $\theta \in [-k,0]$ and $i\ge n$.

\textbf{\textbf{\emph{Step}}\,\emph{1:}} Firstly, we will give a priori estimate for $u$. Multiplying $(\ref{3.3-3})_{1}$ by ${\chi _{i,{\rm{ }}j}}(t)$ and summing $j$ from $1$ to $i$, we deduce
\begin{equation}
\begin{aligned}
&\frac{d}{d t}(\left\|u_{i}(t)\right\|^{2}+\varepsilon(t)\left\|\nabla u_{i}(t)\right\|^{2})+(2a(l(u_{i}))-\varepsilon^{\prime}(t))\left\|\nabla u_{i}(t)\right\|^{2}\\
&= 2\left(f\left(u_{i}(t)\right)+g(t, u_{t,i})+h(t), u_{i}(t)\right).
\end{aligned}
\label{3.4-3}
\end{equation}

Using (\ref{1.6-3}), we obtain there exists constants $0<\delta_1<\frac{m \lambda_1}{2}$ and $C_{\delta_1}>0$ such that
\begin{equation}
2\left(f\left(u_i(t)\right), u_i(t)\right) \leqslant\left(m \lambda_1-2 \delta_1\right)\left\|u_i(t)\right\|^2+2C_{\delta_1}.
\label{3.5-3}
\end{equation}

From $(\ref{1.11-3})-(\ref{1.13-3})$ and the Young inequality, we conclude
\begin{equation}
2\left(g\left(t, u_{t, i}\right), u_i(t)\right) \leqslant \frac{2C_g}{m \lambda_1}\left\|u_{t, i}\right\|_{C_{L^2(\Omega)}}^2+\frac{m \lambda_1}{2}\left\|u_i(t)\right\|^2.
\label{3.6-3}
\end{equation}

Similarly, by the Young inequality, we derive
\begin{equation}
2\left(h(t), u_i(t)\right) \leqslant \frac{2}{m \lambda_1}\|h(t)\|^2+\frac{m \lambda_1}{2}\left\|u_i(t)\right\|^2.
\label{3.7-3}
\end{equation}

Then substituting (\ref{3.5-3})$-$(\ref{3.7-3}) into (\ref{3.4-3}), it follows that
\begin{equation}
\begin{aligned}
&\frac{d}{d t}(\|u_{i}(t)\|^{2}+\varepsilon(t)\|\nabla u_{i}(t)\|^{2})+(2 a(l(u_{i}))-\varepsilon^{\prime}(t))\|\nabla u_{i}(t)\|^{2} \\
& \leqslant\left(2 m \lambda_1-2 \delta_1\right)\left\|u_i(t)\right\|^2+\frac{2 C_g}{m \lambda_1}\left\|u_{t,i}\right\|_{C_{L^2(\Omega)}}^2+\frac{2}{m \lambda_1}\|h(t)\|^2+2C_{\delta_1}.
\end{aligned}
\label{3.8-3}
\end{equation}

Thanks to the Poincar\'{e} inequality, we arrive at
\begin{equation}
\begin{aligned}
&\frac{d}{d t}(\|u_{i}(t)\|^{2}+\varepsilon(t)\|\nabla u_{i}(t)\|^{2})+(2 a(l(u_{i}))-2m-\varepsilon^{\prime}(t))\|\nabla u_{i}(t)\|^{2}+2\delta_1\|u_{i}(t)\|^{2} \\
& \leqslant\frac{2 C_g}{m \lambda_1}\left\|u_{t,i}\right\|_{C_{L^2(\Omega)}}^2+\frac{2}{m \lambda_1}\|h(t)\|^2+2C_{\delta_1}.
\end{aligned}
\label{3.9-3}
\end{equation}

When $\varepsilon(t)$ is decreasing, from $(\ref{1.2-3})$ and $(\ref{1.3-3})$, we obtain $-\varepsilon^{\prime}(t)\in [0, L]$, then combine it with (\ref{1.4-3}) yields
\begin{equation}
\begin{aligned}
& \min \left.\{\left.(2 a(l(u_{i}))-2m-\varepsilon^{\prime}(t))\|\nabla u_{i}(t)\|^{2}\right.\}\right.=0
\end{aligned}
\label{3.10-3}
\end{equation}
and
\begin{equation}
\max \left.\{\left.(2 a(l(u_{i}))-2m-\varepsilon^{\prime}(t))\|\nabla u_{i}(t)\|^{2}\right.\}\right. =(2( M-m)+L)\left\|\nabla u_i(t)\right\|^2>0.
\label{3.11-3}
\end{equation}

Besides, when $\varepsilon(t)$ is increasing, from $(\ref{1.2-3})$ and $(\ref{1.3-3})$, we conclude $-\varepsilon^{\prime}(t)\in [-L, 0]$, then by (\ref{1.5-3}), we derive
\begin{equation}
\min \left.\{\left.(2 a(l(u_{i}))-2m-\varepsilon^{\prime}(t))\|\nabla u_{i}(t)\|^{2}\right.\}\right. =L\left\|\nabla u_i(t)\right\|^2>0
\label{3.12-3}
\end{equation}
and
\begin{equation}
\max \left.\{\left.(2 a(l(u_{i}))-2m-\varepsilon^{\prime}(t))\|\nabla u_{i}(t)\|^{2}\right.\}\right. =2(M-m)\left\|\nabla u_i(t)\right\|^2>0.
\label{3.13-3}
\end{equation}

Let $\widetilde m=2(M-m)$ and by $(\ref{3.9-3})-(\ref{3.13-3})$, we obtain
\begin{equation}
\begin{aligned}
&\frac{d}{d t}(\|u_{i}(t)\|^{2}+\varepsilon(t)\|\nabla u_{i}(t)\|^{2})+\widetilde m \|\nabla u_{i}(t)\|^{2}+2\delta_1\|u_{i}(t)\|^{2} \\
& \leqslant\frac{2 C_g}{m \lambda_1}\left\|u_{t,i}\right\|_{C_{L^2(\Omega)}}^2+\frac{2}{m \lambda_1}\|h(t)\|^2+2C_{\delta_1}.
\end{aligned}
\label{3.14-3}
\end{equation}

Integrating (\ref{3.14-3}) from $\tau$ to $t$ and using $(\ref{1.2-3})-(\ref{1.5-3})$, we deduce
\begin{equation}
\begin{aligned}
& \left.\|u_i(t)\right\|^2+|\varepsilon(t)|\| \nabla u_i(t)\|^2+\widetilde{m} \int_\tau^t\| \nabla u_i(s)\|^2 d s+2 \delta_1 \int_\tau^t\| u_i(s) \|^2 d s \\
& \leqslant\left\|u_i(\tau)\right\|^2+\varepsilon(\tau)\left\|\nabla u_i(\tau)\right\|^2+\frac{2 C_g}{m \lambda_1} \int_\tau^t\left\|u_{s, i}\right\|_{C_{L^2(\Omega)}}^2ds+\frac{2}{m \lambda_1} \int_\tau^t\|h(s)\|^2 d s +2 C_{\delta_1}(t-\tau).
\end{aligned}
\label{3.15-3}
\end{equation}

Putting $t+\theta$ instead of $t$ with $\theta \in [-k,0]$ in (\ref{3.15-3}), we arrive at
\begin{equation}
\begin{aligned}
&\left\|u_{t, i}\right\|_{C_{L^{2}(\Omega)}}^{2}+|\varepsilon_{t}|\left\|\nabla u_{t, i}\right\|_{C_{L^{2}(\Omega)}}^{2}+\widetilde{m} \int_{\tau}^{t}\left\|\nabla u_{s, i}\right\|_{C_{L^{2}(\Omega)}}^{2} d s+2\delta_1 \int_{\tau}^{t}\left\|u_{s, i}\right\|_{C_{L^{2}(\Omega)}}^{2}  d s\\
&\leq\|\phi\|_{C_{L^{2}(\Omega)}}^{2}+\varepsilon_{\tau}\|\nabla \phi\|_{C_{L^{2}(\Omega)}}^{2}+\frac{2 C_g}{m \lambda_1} \int_\tau^t\left\|u_{s, i}\right\|_{C_{L^2(\Omega)}}^2ds+\frac{2}{m \lambda_1} \int_\tau^t\|h(s)\|^2 d s +2 C_{\delta_1}(t-\tau).
\end{aligned}
\label{3.16-3}
\end{equation}

Then from the Gronwall inequality, it follows that
\begin{align}
\left\|u_{t, i}\right\|_{C_{L^{2}(\Omega)}}^{2}+|\varepsilon_{t}|\left\|\nabla u_{t, i}\right\|_{C_{L^{2}(\Omega)}}^{2}&\leq e^{\frac{2C_g}{m \lambda_1}(t-\tau)} \left(\|\phi\|_{C_{L^{2}(\Omega)}}^{2}+\varepsilon_{\tau}\|\nabla \phi\|_{C_{L^{2}(\Omega)}}^{2}\right)\non\\
& +\frac{2}{m \lambda_1}e^{\frac{2C_g}{m \lambda_1}(t-\tau)} \int_\tau^t\|h(s)\|^2 d s \non\\
&+2e^{\frac{2C_g}{m \lambda_1}(t-\tau)}  C_{\delta_1}(t-\tau).
\label{3.17-3}
\end{align}

Thanks to (\ref{3.16-3}) and (\ref{3.17-3}), we conclude
\begin{equation}
\left\{u_{i}\right\} \text { is bounded in } C([\tau-k, T]; \mathcal{H}_{t}(\Omega)) \cap L^{2}(\tau, T ; {H}_{0}^{1}(\Omega)).
\label{3.18-3}
\end{equation}

Furthermore, by $(\ref{1.6-3})-(\ref{1.10-3})$, we obtain
\begin{equation}
f(u_{i}(t)) \text { is bounded in } L^{q}(\tau, T ; L^{q}(\Omega)),
\label{3.19-3}
\end{equation}
where $q=\frac{p}{p-1}$ with $p \ge 2$.

Multiplying $(\ref{3.3-3})_{1}$ by $\chi_{i, j}^{\prime}(t)$ and summing $j$ from 1 to $i$, using $(\ref{1.2-3})-(\ref{1.5-3})$ and doing some similar calculations to $(\ref{3.10-3})-(\ref{3.13-3})$,  we derive
\begin{equation}
\begin{aligned}
2\|u_{i}^{\prime}(t)\|^{2}+2|\varepsilon(t)|\|\nabla u_{i}^{\prime}(t)\|^{2}+\frac{M}{2} \frac{d}{d t}\left\|\nabla u_{i}(t)\right\|^{2} \leq 2(f(u_{i}(t))+g(t, u_{t, i})+h(t), u_{i}^{\prime}(t)).
\end{aligned}
\label{3.20-3}
\end{equation}

Using the Young inequality, we arrive at
\begin{equation}
2(f(u_{i}(t)), u_{i}^{\prime}(t))\leq 2\|f(u_{i}(t))\|^{2}+\frac{1}{2}\|\nabla u_{i}^{\prime}(t)\|^{2}
\label{3.21-3}
\end{equation}
and
\begin{equation}
2(h(t), u_{i}^{\prime}(t)) \leq 2\|h(t)\|^{2}+\frac{1}{2}\|\nabla u_{i}^{\prime}(t)\|^{2}.
\label{3.22-3}
\end{equation}

In addition, by $(\ref{1.2-3})$, $(\ref{1.3-3})$ and the Young inequality, we conclude
\begin{equation}
2(g(t, u_{t, i}), u_{i}^{\prime}(t)) \leq 2C_{g}\|u_{t, i}\|_{C_{L^{2}(\Omega)}}^{2}+\frac{1}{2}\|u_{i}^{\prime}(t)\|^{2}.
\label{3.23-3}
\end{equation}

Substituting $(\ref{3.21-3})-(\ref{3.23-3})$ into (\ref{3.20-3}), we obtain
\begin{equation}
\frac{1}{2}\|u_{i}^{\prime}(t)\|^{2}+2|\varepsilon(t)|\|\nabla u_{i}^{\prime}(t)\|^{2}+\frac{M}{2} \frac{d}{d t}\|\nabla u_{i}(t)\|^{2}
\leq 2\|f(u_{i}(t))\|^{2}+2C_{g}\|u_{t, i}\|_{C_{L^{2}(\Omega)}}^{2}+2\|h(t)\|^{2}.
\label{3.24-3}
\end{equation}

Integrating (\ref{3.24-3}) from $\tau$ to $t$, we deduce
\begin{equation}
\begin{aligned}
&\frac{1}{2}\|u_{i}(t)\|^{2}+2\int_{\tau}^{t} |\varepsilon(s)|\left\|\nabla u_{i}^{\prime}(s)\right\|^{2} d s+\frac{M}{2} \int_{\tau}^{t}\|\nabla u_{i}^{\prime}(s)\|^{2} d s\\
&\leq\frac{1}{2}\|u_{i}(\tau)\|^{2}+2 \int_{\tau}^{t} \|f(u_{i}(t))\|^{2} d s+2C_{g} \int_{\tau}^{t}\|u_{s, i}\|_{C_{L^{2}(\Omega)}}^{2} d s+2\int_{\tau}^{t}\|h(s)\|^{2} d s.
\end{aligned}
\label{3.25-3}
\end{equation}

Putting $t+\theta$ instead of $t$ with $\theta \in [-k, 0]$ in $(\ref{3.25-3})$, using the Gronwall inequality and by similar calculations to $(\ref{3.16-3})$ and (\ref{3.17-3}), we arrive at
\begin{equation}
\left\{u_{i}\right\} \text { is bounded in } L^{\infty}(\tau, T ; \mathcal H_{t}(\Omega)) \cap L^{2}(\tau,T ; \mathcal H_{0}^{1}(\Omega))
\label{3.26-3}
\end{equation}
and
\begin{equation}
\{\partial_{t} u_{i}\} \text { is bounded in } L^{2}(\tau,T ; \mathcal H_{t}(\Omega)).
\label{3.27-3}
\end{equation}

By (\ref{3.18-3}), (\ref{3.19-3}), (\ref{3.26-3}), (\ref{3.27-3}) and the Aubin-Lions lemma, we derive that there exists a subsequence $\left\{\tilde u_{i}\right\}$ of $\left\{u_{i}\right\}$, $\tilde u \in L^{\infty}\left(\tau, T ; \mathcal{H}_{t}(\Omega)\right) \cap L^{2}(\tau, T ; H_{0}^{1}(\Omega))$ and $\partial_{t} \tilde u \in L^{2}\left(\tau, T ; \mathcal{H}_{t}(\Omega)\right)$ such that
\begin{equation}\label{3.28-3}
\tilde u_{i} \rightharpoonup u \quad \text { weakly-star in } L^{\infty}(\tau-k, T ;\mathcal H_{t}(\Omega)),
\end{equation}
\begin{equation}
\tilde u_{i} \rightharpoonup u \quad \text { weakly in } L^{2}(\tau, T ; H_0^1(\Omega)),
\label{3.29-3}
\end{equation}
\begin{equation}
\tilde u_{i} \rightharpoonup u \quad \text { weakly in } L^{p}(\tau, T ; L^{p}(\Omega)),
\label{3.30-3}
\end{equation}
\begin{equation}
f(\tilde u_{i}) \rightharpoonup f(u) \quad \text { weakly in } L^{q}(\tau, T ; L^{q}(\Omega)),
\label{3.31-3}
\end{equation}
\begin{equation}
a\left(l(\tilde u_{i})\right) \tilde u_{i} \rightharpoonup a\left(l(u)\right) u \quad \text { weakly in } L^{2}(\tau, T ; \mathcal H_{t}(\Omega)),
\label{3.32-3}
\end{equation}
\begin{equation}
{\partial _t}{\tilde u_i} \rightharpoonup {\partial _t}u \quad \text { weakly in } L^{2}(\tau, T ;\mathcal H_{t}(\Omega)),
\label{3.33-3}
\end{equation}
\begin{equation}
\tilde u_{i} \rightarrow u \quad \text{ in } C([\tau-k, T];\mathcal H_{t}(\Omega)).
\label{3.34-3}
\end{equation}

\textbf{\textbf{\emph{Step}}\,\emph{2:}} By $(\ref{3.28-3})-(\ref{3.34-3})$, noticing $\left\{e_{j}\right\}_{j=1}^n$ is dense in $H_{0}^{1}(\Omega) \cap L^{p}(\Omega)$ and letting $\zeta \in C^{1}([\tau, T] ; H_{0}^{1}(\Omega))$ with $\zeta(T)=0$ be the test function, then we can easily derive the initial value of $u$ satisfies $u_{\tau, i}=\phi$.

From the above steps, we obtain that $u$ is a weak solution to problem $(\ref{1.1-3})$.

\textbf{\textbf{\emph{Step}}\,\emph{3:}}
Finally, we will verify the continuity of $u$. Suppose $w_{1}$ and $w_{2}$ are solutions to problem $(\ref{1.1-3})$ and their initial values are $w_{1}(\tau+\theta)$ and $w_{2}(\tau+\theta)$ respectively, then we derive
\begin{equation}
\left\{\begin{array}{ll}
\partial_{t}w_{1}-\varepsilon(t) \Delta \partial_{t}w_{1}-a(l(w_{1})) \Delta w_{1}=f(w_{1})+g(t,w_{t,1})+h(t) & \text { in } \Omega \times(\tau, \infty), \\
w_{1}(x,t)=0 & \text { on } \partial \Omega\times(\tau, \infty), \\
w_{1}(x, \tau+\theta)=\phi_{1}(x,\theta),  &\,\, x \in \Omega,\, \theta\in[-k,0],
\end{array}\right.
\label{3.35-3}
\end{equation}
and
\begin{equation}
\left\{\begin{array}{ll}
\partial_{t}w_{2}-\varepsilon(t) \Delta \partial_{t}w_{2}-a(l(w_{2})) \Delta w_{2}=f(w_{2})+g(t, w_{t,2})+h(t) & \text { in } \Omega \times(\tau, \infty), \\
w_{2}(x,t)=0 & \text { on } \partial \Omega\times(\tau, \infty), \\
w_{2}(x, \tau+\theta)=\phi_{2}(x,\theta),  &\,\, x \in \Omega,\,\theta\in[-k,0].
\end{array}\right.
\label{3.36-3}
\end{equation}

Subtracting $(\ref{3.36-3})_{1}$ from $(\ref{3.35-3})_{1}$ and letting $u=w_{1}-w_{2}$, then taking $L^2(\Omega)$-inner product between $u$ and the resulting equation, we conclude
\begin{equation}
\begin{aligned}
&\frac{d}{d t}(\|u\|^{2}+\varepsilon(t)\|\nabla u\|^{2})+2(a(l(u_{1}))-\varepsilon^{\prime}(t))\|\nabla u\|^{2} \\
&=2(a(l(w_{2}))-a(l(w_{1})))(\nabla w_{2}, \nabla u)+2(f(w_{1})-f(w_{2}), u)+2(g(t,w_{t,1})-g(t,w_{t,2}),u).
\end{aligned}
\label{3.37-3}
\end{equation}

Thanks to the Poincar\'{e} inequality and doing similar calculations to those in $(\ref{3.5-3})-(\ref{3.14-3})$, we deduce
\begin{equation}
\frac{d}{d t}(\|u\|^{2}+|\varepsilon(t)|\|\nabla u \|^{2}) \leq C\left(\|u_{t}\||_{C_{L^{2}(\Omega)}}^{2}+|\varepsilon(t)|\|\nabla u_{t}\|_{C_{L^{2}(\Omega)}}^{2}\right).
\label{3.38-3}
\end{equation}

Integrating $(\ref{3.38-3})$ from $\tau$ and $t$ and putting $t+\theta$ instead of $t$ with $\theta \in [-k,0]$ in the obtained equation, we arrive at
\begin{equation}
\begin{aligned}
\|u_{t}\|_{C_{L^{2}(\Omega)}}^{2}+|\varepsilon_{t}|\|\nabla u_{t}\|_{C_{L^{2}(\Omega)}}^{2} &\leq \|u_{\tau}\|_{C_{L^{2}(\Omega)}}^{2}+|\varepsilon_{\tau}|\|\nabla u_{\tau}\|_{C_{L^{2}(\Omega)}}^{2}\\
&+C\int_{\tau}^{t}\left(\left\|u_{s}\right\|_{C_{L^{2}(\Omega)}}^{2}+|\varepsilon_{s}|\|\nabla u_{s}\|_{C_{L^{2}(\Omega)}}^{2}\right) d s.
\end{aligned}
\label{3.39-3}
\end{equation}

From the Gronwall inequality, we obtain
\begin{equation}
\|u_{t}\|_{C_{L^{2}(\Omega)}}^{2}+|\varepsilon_{t}|\|\nabla u_{t}\|_{C_{L^{2}(\Omega)}}^{2} \leq e^{C(t-\tau)}(\|u_{\tau}\|_{C_{L^{2}(\Omega)}}^{2}+\varepsilon_{\tau}\|\nabla u_{\tau}\|_{C_{L^{2}(\Omega)}}^{2}).
\label{3.40-3}
\end{equation}

Therefore, it follows the uniqueness and continuity of the weak solution $u$ to problem $(\ref{1.1-3})$. $\hfill$$\Box$

As a result, by Theorem $\rm{\ref{th3.3-3}}$, we derive problem $(\ref{1.1-3})$ generates a continuous process $
U(t, \tau):C_{\mathcal{ H}_{\tau}(\Omega)} \rightarrow C_{\mathcal{ H}_{t}(\Omega)}$ and $U(t, \tau)\phi=u(t)$ is a unique weak solution to problem $(\ref{1.1-3})$.

\section{Existence and preliminary estimates of regularity for pullback attractors}
In this section, based on Definition \ref{def2.14-3}, we will use Lemmas \ref{lem2.9-3} and \ref{lem2.15-3} to investigate the existence of pullback $\mathcal D_{C_{\mathcal{H}_t(\Omega)}}$-attractors in the time-dependent space $C_{\mathcal{H}_t(\Omega)}$.

Firstly, we shall prove the following lemma, which is beneficial to discuss the existence of a pullback $\mathcal D_{C_{\mathcal{H}_t(\Omega)}}$-absorbing set.
\begin{Lemma} \label{lem4.1-3}
Under the assumptions of $\varepsilon(t)$, $a(l(u))$, $f$, $g$ and $h$ in $\S 1$, assume $\phi \in C_{\mathcal H_{t}(\Omega)}$ is given, if the parameter $m$ in hypothesis $(\ref{1.4-3})$ and $(\ref{1.5-3})$ further satisfies $m>\frac{3}{2}+\frac{L}{2}+\frac{1}{4 \lambda_1}$, then the weak solution $u$ to problem $(\ref{1.1-3})$ satisfies
\begin{equation}
\left\|u_{t}\right\|_{C_{L^{2}(\Omega)}}^{2}+\left|\varepsilon_{t}|\|\nabla u_{t}\right\|_{C_{L^{2}(\Omega)}}^{2} \leq R^2_0(t)
\label{4.1-3}
\end{equation}
for any $\tau \leq t \in \mathbb R$, where
\begin{equation}
\begin{aligned}
 R_0^2(t)&=\left(1+\frac{2C_g e^{\beta k}}{\left(1+\lambda_1L\right)\left(\beta-\beta_1\right)}\right)\left(\|\phi\|_{C_{L^{2}(\Omega)}}^{2}+\varepsilon_{\tau}\|\nabla \phi\|^{2}_{C_{L^{2}(\Omega)}}\right) e^{-\beta_{1}(t-\tau)} \\
& +\frac{1}{\delta}\left(1+\frac{2 C_g e^{\beta k}}{\left(1+\lambda_1 L\right)\left(\beta-\beta_1\right)}\right) e^{\beta k} \int_\tau^t e^{-\beta_1(t-s)}\|h(s)\|^2 d s,
\end{aligned}
\label{4.2-3}
\end{equation}
$\beta_1=\beta-\frac{2 C_g}{1+\lambda_1 L} e^{\beta k}>0$, $0<\beta \leqslant \frac{\lambda_1 \bar{\delta}}{1+\lambda_{1}L}$ when $\varepsilon(t)$ is decreasing, $0<\beta \leqslant \frac{\bar{\delta}}{\lambda_1^{-1}+\alpha}$ when $\varepsilon(t)$ is increasing, $\alpha$ is the same as in $(\ref{1.2-3})$ and $\bar{\delta}>L$.
\end{Lemma}
$\mathbf{Proof.}$
Taking the inner product of $L^2(\Omega)$ with problem $(\ref{1.1-3})_1$ and $u$, we obtain
\begin{equation}
\begin{aligned}
&\frac{d}{d t}(\|u(t)\|^{2}+\varepsilon(t)\|\nabla u(t)\|^{2})+(2 a(l(u))-\varepsilon^{\prime}(t))\|\nabla u\|^{2}=2(f(u)+g(t,u_{t})+h(t), u).
\label{4.3-3}
\end{aligned}
\end{equation}

Using (\ref{1.6-3}), we deduce
\begin{equation}
2(f(u), u) \leq 2\lambda_1\|u\|^2.
\label{4.4-3}
\end{equation}

By the Young inequality and $(\ref{1.11-3})-(\ref{1.13-3})$, we conclude
\begin{equation}
2(g(t, u_{t}), u) \leq 2 C_{g} \|u_{t}\|^{2}_{C_{L^{2}(\Omega)}}+\frac{1}{2}\|u\|^{2}.
\label{4.5-3}
\end{equation}

Similarly, from the Young inequality, it follows that there exists a constant $0<\delta<\lambda_1$ such that
\begin{equation}
2(h(t), u) \leq \frac{1}{\delta}\|h(t)\|^{2}+\delta\| u\|^{2}.
\label{4.6-3}
\end{equation}

Then inserting $(\ref{4.4-3})-(\ref{4.6-3})$ into (\ref{4.3-3}) and using the Poincar\'{e} inequality yield
\begin{equation}
\begin{aligned}
&\frac{d}{d t}\left(\|u\|^{2}+\varepsilon(t)\|\nabla u\|^{2}\right)+(2 a(l(u))-\varepsilon^{\prime}(t)-2-\frac{1}{2 \lambda_1}-\frac{\delta}{\lambda_1})\|\nabla u\|^{2}\\
&\leq 2C_{g}\|u_{t}\|_{C{_{L^{2}(\Omega)}}} ^{2}+\frac{1}{\delta}\|h(t)\|^{2}.
\end{aligned}
\label{4.7-3}
\end{equation}

When $\varepsilon(t)$ is decreasing, noticing $-\varepsilon^{\prime}(t)\in [0, L]$, then by $m>\frac{3}{2}+\frac{L}{2}+\frac{1}{4 \lambda_1}$ and $0<\delta<\lambda_1$, we derive
\begin{equation}
\min \left.\{\left.(2 a(l(u))-\varepsilon^{\prime}(t)-2-\frac{1}{2\lambda_1}-\frac{\delta}{\lambda_1})\|\nabla u(t)\|^{2}\right.\}\right. =L\left\|\nabla u(t)\right\|^2>0.
\label{4.8-3}
\end{equation}

Besides, when $\varepsilon(t)$ is increasing and noticing $-\varepsilon^{\prime}(t)\in [-L, 0]$, we obtain
\begin{equation}
\min \left.\{\left.(2 a(l(u))-\varepsilon^{\prime}(t)-2-\frac{1}{2\lambda_1}-\frac{\delta}{\lambda_1})\|\nabla u(t)\|^{2}\right.\}\right. =0.
\label{4.9-3}
\end{equation}

From $(\ref{4.7-3})-(\ref{4.9-3})$, it follows that there exists a constant $\bar{\delta}>L$ such that
\begin{equation}
\begin{aligned}
&\frac{d}{d t}\left(\|u\|^{2}+\varepsilon(t)\|\nabla u\|^{2}\right)+\bar{\delta}\|\nabla u\|^{2}
\leq 2C_{g}\|u_{t}\|_{C{_{L^{2}(\Omega)}}} ^{2}+\frac{1}{\delta}\|h(t)\|^{2}.
\end{aligned}
\label{4.10-3}
\end{equation}

Using $(\ref{1.2-3})-(\ref{1.4-3})$ and the Poincar\'{e} inequality, then
there exists a constant $\beta$ satisfies $0<\beta \leqslant \frac{\lambda_1 \bar{\delta}}{1+\lambda_{1}L}$ when $\varepsilon(t)$ is decreasing and $0<\beta \leqslant \frac{\bar{\delta}}{\lambda_1^{-1}+\alpha} $ when $\varepsilon(t)$ is increasing such that
\begin{equation}
\begin{aligned}
&\frac{d}{d t}\left(\|u\|^{2}+|\varepsilon(t)|\|\nabla u\|^{2}\right)+\beta(\|u\|^{2}+|\varepsilon(t)|\|\nabla u\|^{2})
\leq 2C_{g}\|u_{t}\|_{C{_{L^{2}(\Omega)}}} ^{2}+\frac{1}{\delta}\|h(t)\|^{2}.
\end{aligned}
\label{4.11-3}
\end{equation}

Furthermore, multiplying (\ref{4.11-3}) by $e^{\beta t}$ and integrating the resulting inequality from $\tau$ to $t$, we deduce
\begin{equation}
\begin{aligned}
e^{\beta t}\left(\|u\|^{2}+|\varepsilon(t)|\|\nabla u\|^{2}\right) & \leq e^{\beta \tau}\left(\|u(\tau)\|^{2}+\varepsilon(\tau)\|\nabla u(\tau)\|^{2}\right)\\
&+2C_{g} \int_{\tau}^{t} e^{\beta s}\|u_{s}\|_{C_{L^{2}(\Omega)}}^{2} d s+\frac{1}{\delta}\int_{\tau}^{t} e^{\beta s}\|h(s)\|^{2} d s.
\end{aligned}
\label{4.12-3}
\end{equation}

Putting $t+\theta$ instead of $t$ with $\theta \in [-k,0]$ in (\ref{4.12-3}), then noticing $u_t=u(t+\theta)$, $u(x, \tau+\theta)=\phi(x, \theta)$ and $\varepsilon(\tau+\theta)=\varepsilon_\tau$, we arrive at
\begin{equation}
\begin{aligned}
e^{\beta(t+\theta)}\left(\|u_{t}\|_{C_{L^{2}(\Omega)}}^{2}+|\varepsilon_{t}|\|\nabla u_{t}\|^{2}_{C_{L^{2}(\Omega)}}\right) &\leq \left(\|\phi\|_{C_{L^{2}(\Omega)}}^{2}+\varepsilon_{\tau}\|\nabla \phi\|_{C_{L^{2}(\Omega)}}^{2}\right)e^{\beta(\tau+\theta)}\\
&+2C_{g} \int_{\tau}^{t} e^{\beta s}\|u_{s}\|_{C_{L^{2}(\Omega)}}^{2} d s+\frac{1}{\delta}\int_{\tau}^{t} e^{\beta s}\|h(s)\|^{2} d s.
\end{aligned}
\label{4.13-3}
\end{equation}

Taking $\theta=-k$ in (\ref{4.13-3}), we obtain
\begin{equation}
\begin{aligned}
e^{\beta t}\left(\|u_{t}\|_{C_{L^{2}(\Omega)}}^{2}+|\varepsilon_{t}|\|\nabla u_{t}\|^{2}_{C_{L^{2}(\Omega)}}\right) &\leq \left(\|\phi\|_{C_{L^{2}(\Omega)}}^{2}+\varepsilon_{\tau}\|\nabla \phi\|_{C_{L^{2}(\Omega)}}^{2}\right)e^{\beta\tau}\\
&+2C_{g}e^{\beta k} \int_{\tau}^{t} e^{\beta s}\|u_{s}\|_{C_{L^{2}(\Omega)}}^{2} d s+\frac{1}{\delta}e^{\beta k}\int_{\tau}^{t} e^{\beta s}\|h(s)\|^{2} d s.
\end{aligned}
\label{4.14-3}
\end{equation}

Then from $(\ref{1.2-3})-(\ref{1.3-3})$ and the Poincar\'{e} inequality, we conclude
\begin{equation}
\begin{aligned}
\|u_{t}\|_{C_{L^{2}(\Omega)}}^{2}e^{\beta t} & \leq \frac{1}{1+\lambda_{1}L} \left(\|\phi\|_{C_{L^{2}(\Omega)}}^{2}+\varepsilon_{\tau}\|\nabla \phi\|_{C_{L^{2}(\Omega)}}^{2}\right)e^{\beta \tau}\\
&+\frac{2C_{g}}{1+\lambda_{1}L} e^{\beta k} \int_{\tau}^{t} e^{\beta s}\|u_{s}\|_{C_{L^{2}(\Omega)}}^{2} d s+ \frac{1}{(1+\lambda_{1}L)\delta}e^{\beta k}\int_{\tau}^{t} e^{\beta s}\|h(s)\|^{2} d s.
\end{aligned}
\label{4.15-3}
\end{equation}

Next, we shall use Lemma \ref{lem2.16-3} to prove $u$ satisfies $(\ref{4.1-3})$ and $(\ref{4.2-3})$.

Suppose
\begin{equation}
\bar w(t)=\left\|u_{t}\right\|_{C_{L^{2}(\Omega)}}^{2}e^{\beta t},
\label{4.16-3}
\end{equation}
\begin{equation}
\bar{u}(s)=\frac{2C_g}{1+\lambda_{1}L} e^{\beta k}
\label{4.17-3}
\end{equation}
and
\begin{equation}
\begin{aligned}
\bar v(t) & = \frac{1}{1+\lambda_{1}L}\left(\|\phi\|_{C_{L^{2}(\Omega)}}^{2}+\varepsilon_{\tau}\|\nabla \phi\|_{C_{L^{2}(\Omega)}}^{2}\right) e^{\beta \tau}+ \frac{1}{(1+\lambda_{1}L)\delta}e^{\beta k}\int_{\tau}^{t} e^{\beta s}\|h(s)\|^{2} d s.
\end{aligned}
\label{4.18-3}
\end{equation}

Let $a=\tau$ and from $(\ref{4.17-3})$ and $(\ref{4.18-3})$, we derive
\begin{equation}
\bar v(a) e^{\int_{a}^{t}\bar u(s) d s}=\frac{1}{1+\lambda_{1} L} \left(\|\phi\|_{C_{L^{2}(\Omega)}}^{2}+\varepsilon_{\tau}\|\nabla \phi\|_{C_{L^{2}(\Omega)}}^{2}\right)e^{\beta \tau} e^{\frac{2C_{g}}{1+\lambda_{1} L} e^{\beta k} \cdot(t-\tau)}.
\label{4.19-3}
\end{equation}

By (\ref{4.18-3}), it follows
\begin{equation}
\frac{d \bar{v}(s)}{d s}=\frac{1}{\left(1+\lambda_1 L\right) \delta} e^{\beta(s+k)}\|h(s)\|^2.
\label{4.20-3}
\end{equation}

Thanks to (\ref{4.19-3}) and (\ref{4.20-3}), we deduce
\begin{equation}
\int_{a}^{t}e^{\int_{s}^{t} \bar{u}(r)d r} \cdot \frac{d \bar{v}(s)}{d s} d s=\frac{1}{\left(1+\lambda_1 L)\delta\right.} e^{\beta k} e^{\frac{2C_g}{1+\lambda_1 L} e^{\beta k} \cdot t} \int_\tau^t e^{s\left(\beta-\frac{2C_g}{1+\lambda_1 L} e^{\beta k}\right)}\|h(s)\|^2 d s.
\label{4.21-3}
\end{equation}

Let $\beta_1=\beta-\frac{2 C_g}{1+\lambda_1 L} e^{\beta k}>0$, then by $(\ref{4.19-3})$ and $(\ref{4.21-3})$ and Lemma \ref{lem2.16-3}, we obtain
\begin{align}
\bar w(t)=\|u_{t}\|_{C_{L^{2}(\Omega)}}^{2}e^{\beta t} & \leq \bar v(a) e^{\int_{a}^{t}\bar u(s) d s}+\int_{a}^{t} e^{\int_{s}^{t} \bar    u(r) d r} \cdot \frac{d \bar v(s)}{d s} d s\non\\
&=\frac{1}{1+\lambda_{1} L} \left(\|\phi\|_{C_{L^{2}(\Omega)}}^{2}+\varepsilon_{\tau}\|\nabla \phi\|_{C_{L^{2}(\Omega)}}^{2}\right)e^{\beta \tau} e^{\frac{2C_{g}}{1+\lambda_{1} L} e^{\beta k} \cdot(t-\tau)}\non\\
&+\frac{1}{\left(1+\lambda_1 L)\delta\right.} e^{\beta k} e^{\frac{2C_g}{1+\lambda_1 L} e^{\beta k} \cdot t} \int_\tau^t e^{\beta_1 s}\|h(s)\|^2 d s.
\label{4.22-3}
\end{align}

Then from (\ref{4.14-3}) and noticing $0<\beta_1<\beta$, we arrive at
\begin{equation}
\|u_{t}\|_{C_{L^{2}(\Omega)}}^{2}+|\varepsilon_{t}|\|\nabla u_{t}\|_{C_{L^{2}(\Omega)}}^{2}\leq (A_1)+(B_1)+2 C_{g}e^{\beta (k-t)} \int_{\tau}^{t} e^{\beta s}\left\|u_{s}\right\|_{C_{L^{2}(\Omega)}}^{2}ds,
\label{4.23-3}
\end{equation}
where
\begin{equation}
(A_1)=\left(\|\phi\|_{C_{L^{2}(\Omega)}}^{2}+\varepsilon_{\tau}\|\nabla \phi\|_{C_{L^{2}(\Omega)}}^{2}\right)e^{-\beta_1(\tau-t)}
\label{4.24-3}
\end{equation}
and
\begin{equation}
(B_1)=\frac{1}{\delta}e^{\beta k} \int_{\tau}^{t} e^{-\beta_{1}(t-s)}\|h(s)\|^{2} d s.
\label{4.25-3}
\end{equation}

Inserting $(\ref{4.22-3})$ into $(\ref{4.23-3})-(\ref{4.25-3})$, we derive
\begin{equation}
\begin{aligned}
\|u_{t}\|_{C_{L^{2}(\Omega)}}^{2}+ & \varepsilon_{t}\|\nabla u_{t}\|_{C_{L^{2}(\Omega)}}^{2}\leq (A_1)+(B_1)+(C_1)+(D_1),
\end{aligned}
\label{4.26-3}
\end{equation}
where
\begin{equation}
(C_1)=2C_{g} e^{\beta(k-t)} \int_{\tau}^{t} \frac{1}{1+\lambda_{1} L} \left(\|\phi\|_{C_{L^{2}(\Omega)}}^{2}+\varepsilon_{\tau}\|\nabla \phi\|_{C_{L^{2}(\Omega)}}^{2}\right)e^{\beta \tau} e^{\frac{2C_{g}}{1+\lambda_{1} L} e^{\beta k} \cdot(s-\tau)} d s
\label{4.27-3}
\end{equation}
and
\begin{equation}
(D_1)=2C_{g} e^{\beta(k-t)} \int_{\tau}^{t} \frac{1}{\delta(1+\lambda_{1} L)} e^{\beta k} e^{\frac{2C_{g}}{1+\lambda_{1} L} e^{\beta k} \cdot s}\left(\int_{\tau}^{s} e^{\beta_{1} s}\|h(s)\|^{2} d s\right) d s.
\label{4.28-3}
\end{equation}

Then by some simple estimates and calculations, we conclude
\begin{equation}
\left(C_1\right) \leq \frac{2 C_g e^{\beta k}}{\left(1+\lambda_1 L\right)\left(\beta-\beta_1\right)}\left(\|\phi\|_{C_{L^2(\Omega)}}^2+\varepsilon_\tau\|\nabla \phi\|_{C_{L^2(\Omega)}}^2\right) e^{-\beta_1(t-\tau)}
\label{4.29-3}
\end{equation}
and
\begin{equation}
\left(D_1\right) \leq \frac{2C_g e^{\beta k}}{\left.\delta(1+\lambda_1 L\right)\left(\beta-\beta_1\right)} e^{\beta k} \int_\tau^t e^{-\beta_1(t-s)}\|h(s)\|^2 d s.
\label{4.30-3}
\end{equation}

By $(\ref{4.26-3})-(\ref{4.27-3})$, we obtain
\begin{equation}
(A_1)+(C_1) \leq\left(1+\frac{2C_{g} e^{\beta k}}{\left(1+\lambda_{1}L)\left(\beta-\beta_{1}\right)\right.}\right)\left(\|\phi\|_{C_{L^{2}(\Omega)}}^{2}+\varepsilon_{\tau}\|\nabla \phi\|_{C_{L^{2}(\Omega)}}^{2}\right) e^{-\beta_{1}(t-\tau)}
\label{4.31-3}
\end{equation}
and
\begin{equation}
(B_1)+(D_1) \leq \frac{1}{\delta}\left(1+\frac{\left. 2C_{g} e^{\beta k}\right.}{\left(1+\lambda_{1} L)\left(\beta-\beta_{1}\right)\right.}\right) e^{\beta k} \int_{\tau}^{t} e^{-\beta_{1}(t-s)}\|h(s)\|^{2}d s.
\label{4.32-3}
\end{equation}

Inserting $(\ref{4.31-3})$ and $(\ref{4.32-3})$ into (\ref{4.26-3}), we conclude $(\ref{4.1-3})$ and $(\ref{4.2-3})$. $\hfill$$\Box$

Next, we shall give the following definition of tempered universe, and then we will give the prove of the existence of a pullback $\mathcal D_{C_{\mathcal{H}_t(\Omega)}}$-absorbing set.
\begin{Definition} \label{def4.2-3}
Suppose $\mathcal D \in \mathcal{H}_t(\Omega)$ is a nonempty set and $\mathcal D_{C_{\mathcal{H}_t(\Omega)}}$ is the class of all families of nonempty sets $\widehat{D}=\{D(t)\}_{t \in \mathbb{R}} \subset \Gamma$ such that
$$
\lim _{\tau \rightarrow-\infty}\left(e^{\beta_{1} \tau} \sup _{u \in D(\tau)}{\|u\|}^{2}_{C_{\mathcal H_{t}(\Omega)}}\right)=0,
$$
where $\Gamma$ is a family of all nonempty subset of the time-dependent space $C_{\mathcal H_{t}(\Omega)}$.
\end{Definition}

\begin{Lemma}\label{lem4.3-3}
Under the assumptions of Lemma $\ref{lem4.1-3}$, if $h$ further satisfies
\begin{equation}
\int_{-\tau}^{t} e^{\beta_{1} s}\|h(s)\|^{2}d s<\infty
\label{4.33-3}
\end{equation}
for any $ \tau < t \in \mathbb{R}$, then the family $\widehat{D}_{0}=\left\{D_{0}(t)\right\}_{t \in \mathbb{R}}$ with $D_{0}(t)=\mathcal{\bar{\mathscr B}}_{C_{\mathcal H_{t}(\Omega)}}\left(0, {R}(t)\right)$, the closed ball in $C_{\mathcal H_{t}(\Omega)}$ of centre zero and radius $R(t)$, where
\begin{equation}
R^2(t)=1+\frac{1}{\delta}\left(1+\frac{2 C_g e^{\beta k}}{\left(1+\lambda_1 L\right)\left(\beta-\beta_1\right)}\right) e^{\beta k} \int_\tau^t e^{-\beta_1(t-s)}\|h(s)\|^2 d s
\label{4.34-3}
\end{equation}
is a pullback $\mathcal D_{C_{\mathcal{H}_t(\Omega)}}$-absorbing set for the process $\{U(t, \tau)\}_{t \geq \tau}$ of weak solutions to problem $(\ref{1.1-3})$ in $C_{\mathcal H_{t}(\Omega)}$. Moreover, $\widehat{D}_0 \in \mathcal D_{C_{\mathcal{H}_t(\Omega)}}$.
\label{lem3.6-3}
\end{Lemma}
$\mathbf{Proof.}$  By Definition \ref{def2.5-3} and Lemma \ref{lem4.1-3}, we conclude $\widehat{D}_{0}$ is a pullback $\mathcal D_{C_{\mathcal{H}_t(\Omega)}}$-absorbing set for the process $\{U(t, \tau)\}_{t \geq \tau}$ of weak solutions to problem $(\ref{1.1-3})$ in $C_{\mathcal H_{t}(\Omega)}$. In addition, from Definition \ref{def4.2-3} and (\ref{4.33-3}), we derive $e^{\beta_{1} t} R^{2}(t) \rightarrow 0$ as $t \rightarrow-\infty$, then $\widehat{D}_0 \in \mathcal D_{C_{\mathcal{H}_t(\Omega)}}$ holds.  $\hfill$$\Box$

Next, we will prove pullback $\mathcal D_{C_{\mathcal{H}_t(\Omega)}}$-$\omega$-limit compact, we further assume $h \in L_b^2\left(\mathbb{R} ; L^2(\Omega)\right)$, then decompose the weak solution $u$ to problem (\ref{1.1-3}) into $u(t)=U(t, \tau) u(\tau)=v_1(t)+v_2(t)$ with $v_1(t)=V_1(t, \tau)v_1(\tau)$ and $v_2(t)=V_2(t, \tau)v_2(\tau)$, which satisfy
\begin{equation}
\left\{\begin{array}{ll}
\partial_{t} v_1-\varepsilon(t)\partial_{t} \Delta v_1-a(l(u)) \Delta v_1=f_0(v_1) & \text { in } \Omega \times(\tau, \infty), \\
v_1(x,t)=0 & \text { on } \partial \Omega \times(\tau, \infty), \\
v_1(x, \tau+s)=\phi(x,s), &\,\, x \in \Omega,\,s\in[-k, 0],
\label{4.35-3}
\end{array}\right.
\end{equation}
and
\begin{equation}
\left\{\begin{array}{ll}
\partial_{t} v_{2}-\varepsilon(t)\partial_{t} \Delta v_{2}-a(l(u)) \Delta v_{2}=f(u)-f_0(v_1)+g(t,u_{t})+h(t) & \text { in } \Omega \times(\tau, \infty), \\
v_{2}(x,t)=0 & \text { on } \partial \Omega \times(\tau, \infty), \\
v_{2}(x, \tau+s)=0, &\,\, x \in \Omega, \,s\in[-k, 0],
\label{4.36-3}
\end{array}\right.
\end{equation}
respectively.

The following lemma will explain why we assume the parameter $\sigma$ satisfies  $0<\sigma<\min \left\{\frac{1}{3}, \frac{n+2-(n-2) \gamma}{2}\right\}$ with $0<\gamma<\frac{n+2}{n-2}$ in $\S 1$, which is crucial to discuss the compactness of the process $\{U(t, \tau)\}_{t \geq \tau}$ in the time-dependent space $C_{\mathcal H_{t}(\Omega)}$.

\begin{Lemma}\label{lem4.4-3}
Assume $0<\sigma<\min \left\{\frac{1}{3}, \frac{n+2-(n-2) \gamma}{2}\right\}$, then it follows the following inequalities
\begin{equation}
\|\cdot\|_{L^2(\Omega)} \leqslant C\|\cdot\|_{L^{\frac{2 n(n-2)}{n(n-4-2\sigma)+4(1+3\sigma)}}} \leqslant C\|\cdot\|_{L^{\frac{2 n}{n-2-2\sigma}}(\Omega)}
\label{4.37-3}
\end{equation}
and
\begin{equation}
\|\cdot\|_{L^{\frac{2n\gamma}{n+2-2\sigma}}(\Omega)} \leqslant C\|\cdot\|_{L^{\frac{2 n}{n-2}}(\Omega)} \leqslant C\| \nabla \cdot \|_{L^2(\Omega)},
\label{4.38-3}
\end{equation}
where $0<\gamma<\frac{n+2}{n-2}$ and $n \ge3$.
\end{Lemma}
$\mathbf{Proof.}$ From the definition of $L^p(\Omega)$, in order to make (\ref{4.37-3}) hold, $\sigma$ should satisfies
\begin{equation}
\frac{2 n(n-2)}{n(n-4-2\sigma)+4(1+3\sigma)}>2
\label{4.39-3}
\end{equation}
and
\begin{equation}
\frac{2 n(n-2)}{n(n-4-2\sigma)+4(1+3\sigma)}<\frac{2 n}{n-2-2\sigma}.
\label{4.40-3}
\end{equation}

Then by (\ref{4.39-3}) and (\ref{4.40-3}), we deduce $\sigma<\frac{1}{3}$.

Similarly, the inequality (\ref{4.38-3}) requires the parameter $\sigma$ to satisfy
\begin{equation}
0<\sigma<\frac{n+2-(n-2) \gamma}{2},
\label{4.41-3}
\end{equation}
where $0<\gamma<\frac{n+2}{n-2}$.

Therefore, by combining $\sigma<\frac{1}{3}$ and  $(\ref{4.41-3})$, we obtain $0<\sigma<\min \left\{\frac{1}{3}, \frac{n+2-(n-2) \gamma}{2}\right\}$ with $0<\gamma<\frac{n+2}{n-2}$ and $n \ge3$.
$\hfill$$\Box$

\begin{Lemma}\label{lem4.5-3}
Under the assumptions of Lemmas $\ref{lem4.3-3}$ and $\ref{lem4.4-3}$, then there exists a constant $\widetilde{C}>0$ depending on $m$, $M$, $l$ and $L$ such that for every $t \geqslant k+\tau$, the weak solution $v_1$ of problem $(\ref{4.35-3})$ satisfies
\begin{equation}
\left\|v_1(t+\theta)\right\|_{C_{\mathcal H_{t}(\Omega)}}^2=\left\|v_{1,t}\right\|_{C_{L^2(\Omega)}}^2+\left|\varepsilon_t\right|\left\|\nabla v_{1, t}\right\|_{C_{L^2(\Omega)}}^2 \leqslant Z\left(\|\phi\|_{C_{\mathcal H_{t}(\Omega)}}^2\right) e^{-\widetilde{C}(t-k-\tau)},
\label{4.42-3}
\end{equation}
where $k>0$ is the same as in $(\ref{1.1-3})_3$ and $Z(\cdot)>0$ is an increasing function.
\end{Lemma}
$\mathbf{Proof.}$ Taking the inner product of $L^2(\Omega)$ with problem $(\ref{4.35-3})_1$ and $v_1$, we conclude
\begin{equation}
\frac{d}{d t}\left(\|v_1\|^2+\varepsilon(t)\|\nabla  v
_1\|^2\right)+\left(2 a(l(u))-\varepsilon^{\prime}(t)\right)\|\nabla  v_1\|^2=(f_0(v_1(t)),v_1(t)).
\label{4.43-3}
\end{equation}

Then by (\ref{1.9-3}) and $u=v_1+v_2$, we arrive at
\begin{equation}
\frac{d}{d t}(\|v_1(t)\|^2+\varepsilon(t)\|\nabla v_1(t)\|^2)+\left(2 a(l (u))-\varepsilon^{\prime}(t)\right)\left\|\nabla v_1(t)\right\|^2+2 l\left\|v_1(t)\right\|^2 \leqslant 0.
\label{4.44-3}
\end{equation}

When $\varepsilon(t)$ is decreasing, from $(\ref{1.2-3})-(\ref{1.4-3})$, we obtain
\begin{equation}
\max \left.\{\left.(2 a(l(u))-\varepsilon^{\prime}(t))\|\nabla v_1\|^{2}\right.\}\right. =(2M+L)\left\|\nabla v_1(t)\right\|^2
\label{4.45-3}
\end{equation}
and
\begin{equation}
\min \left.\{\left.(2 a(l(u))-\varepsilon^{\prime}(t))\|\nabla v_1\|^{2}\right.\}\right. =2m\|\nabla v_1(t)\|^2.
\label{4.46-3}
\end{equation}

Besides, when $\varepsilon(t)$ is increasing, by $(\ref{1.2-3})$, $(\ref{1.3-3})$ and $(\ref{1.5-3})$, we deduce
\begin{equation}
\max \left.\{\left.(2 a(l(u))-\varepsilon^{\prime}(t))\|\nabla v_1\|^{2}\right.\}\right. =2M\|\nabla v_1(t)\|^2
\label{4.47-3}
\end{equation}
and
\begin{equation}
\min \left.\{\left.(2 a(l(u))-\varepsilon^{\prime}(t))\|\nabla v_1\|^{2}\right.\}\right. =(2m+L)\|\nabla v_1(t)\|^2.
\label{4.48-3}
\end{equation}

Inserting $(\ref{4.45-3})-(\ref{4.48-3})$ into $(\ref{4.44-3})$ yields
\begin{equation}
\frac{d}{d t}(\|v_1(t)\|^2+|\varepsilon(t)|\|\nabla v_1(t)\|^2)+(2 M+L)\left\|\nabla v_1(t)\right\|^2+2 l\left\|v_1(t)\right\|^2 \leqslant 0.
\label{4.49-3}
\end{equation}

Then it follows that there exists a constant $\widetilde{C}>0$ depending on $m$, $M$, $l$ and $L$ such that
\begin{equation}
\frac{d}{d t}(\|v_1(t)\|^2+|\varepsilon(t)|\|\nabla v_1(t)\|^2)+\widetilde{C}(\|v_1(t)\|^2+|\varepsilon(t)|\|\nabla v_1(t)\|^2) \leqslant 0.
\label{4.50-3}
\end{equation}

Thanks to the Gronwall inequality, we arrive at
\begin{equation}
\left\|v_1(t)\right\|^2+|\varepsilon(t)|\left\|\nabla v_1(t)\right\|^2 \leqslant e^{-\widetilde{C}(t-\tau)}\left(\left\|v_1(\tau)\right\|^2+\varepsilon(\tau)\left\|\nabla v_1(\tau)\right\|^2\right).
\label{4.51-3}
\end{equation}

Finally, putting $t+\theta$ instead of $t$ with $\theta \in [-k,0]$ in (\ref{4.51-3}), then from $\tau \leqslant t \in \mathbb{R}$ and $k>0$, we obtain there exists an increasing function $Z(\cdot)>0$ such that $(\ref{4.42-3})$ holds. $\hfill$$\Box$

\begin{Lemma}\label{lem4.6-3}
Under the assumptions of Lemmas $\ref{lem4.3-3}$ and $\ref{lem4.4-3}$, assume $\phi \in C_{\mathcal H_t(\Omega)}$ is given, then for every $T_1>0$ there exists a constant $\widetilde{C}_1>0$ depending on $T_1$, $\phi$ and $h$ such that the weak solution $v_2$ of problem $(\ref{4.36-3})$ satisfies
\begin{equation}
\left\|v_2\left(T_1+\tau\right)\right\|_{C_{\mathcal H^1_t(\Omega), \sigma}}^2=\|A^{\frac{\sigma}{2}}v_{2, T_1+\tau}\|_{C_{L^2(\Omega)}}^2+\left|\varepsilon_{T_1+\tau}\right|\|A^{\frac{1+\sigma}{2}} v_{2, T_1+\tau}\|_{C_{L^2(\Omega)}}^2 \leqslant \widetilde{C}_1,
\label{4.52-3}
\end{equation}
where $v_{2, T_1+\tau}=\left(v_2\right)_{T_1+\tau}$, $|\varepsilon_{T_1+\tau}|$ is the absolute value of $\varepsilon(\tau+T_1+\theta)$ with $\theta \in[-k, 0]$ and $\sigma$ is the same as in Lemma $\ref{lem4.4-3}$.
\end{Lemma}
$\mathbf{Proof.}$  Taking the inner product of $L^2(\Omega)$ with problem $(\ref{4.36-3})_1$ and $A^{\sigma}v_2$, then from $f=f_0+f_1$, we obtain
\begin{equation}
\begin{aligned}
& \frac{d}{d t}\left(\|A^{\frac{\sigma}{2}}v_2(t)\|^2+\varepsilon(t)\|A^{\frac{1+\sigma}{2}}v_2(t)\|^2\right)+\left(2 a(l(u))-\varepsilon^{\prime}(t)\right)\|A^{\frac{1+\sigma}{2}}v_2(t)\|^2\\
& \left.=2\left(f(u(t))-f\left(v_1(t)\right), A^{\sigma} v_2(t)\right)+2\left(f_1 (v_1\right)+g(t, u_t)+h(t), A^{\sigma} v_2(t)\right).
\end{aligned}
\label{4.53-3}
\end{equation}

By similar calculations as in Lemma 3.4 in \cite{sy.3}, it is easily to derive
\begin{equation}
\frac{d}{d t}\left(\|A^{\frac{\sigma}{2}}v_2(t)\|^2+\varepsilon(t)\|A^{\frac{1+\sigma}{2}}v_2(t)\|^2\right) \leq C(1+\|A^{\frac{1+\sigma}{2}}v_2(t)\|^2).
\label{4.54-3}
\end{equation}

Then from $(\ref{1.2-3})-(\ref{1.5-3})$ and $(\ref{4.53-3})-(\ref{4.54-3})$, we conclude
\begin{equation}
\frac{d}{d t}\left(\|A^{\frac{\sigma}{2}}v_2(t)\|^2+|\varepsilon(t)|\|A^{\frac{1+\sigma}{2}}v_2(t)\|^2\right)\leq C\left(1+\|A^{\frac{1+\sigma}{2}} v_2(t)\|^2\right)+2(g\left(t, u_t\right)+h(t), A^\sigma v_2(t)).
\label{4.55-3}
\end{equation}

Noticing $0<\sigma<\min \left\{\frac{1}{3}, \frac{n+2-(n-2) \gamma}{2}\right\}$ with $0<\gamma<\frac{n+2}{n-2}$ and by the Young inequality, $(\ref{1.11-3})-(\ref{1.13-3})$, Definition $\ref{def2.10-3}$ and Lemma $\ref{lem2.11-3}$, we deduce
\begin{equation}\label{4.56-3}
2\left(g\left(t, u_t\right), A^{\sigma} v_2(t)\right) \leqslant C_g\left\|u_t\right\|_{C_{L^2(\Omega)}}^2+C \| A^{\frac{1+\sigma}{2} }v_2(t) \|^2
\end{equation}
and
\begin{equation}
2\left(h(t), A^{\sigma} v_2(t)\right) \leqslant\|h(t)\|^2+C\|A^{\frac{1+\sigma}{2}} v_2(t)\|^2.
\label{4.57-3}
\end{equation}

Inserting $(\ref{4.56-3})$ and $(\ref{4.57-3})$ into $(\ref{4.55-3})$, it follows that
\begin{equation}
\frac{d}{d t}\left(\|A^{\frac{\sigma}{2}}v_2(t)\|^2+|\varepsilon(t)|\|A^{\frac{1+\sigma}{2}}v_2(t)\|^2\right) \leqslant C\left(1+\|A^{\frac{1+\sigma}{2}} v_2(t)\|^2\right)+C_g\left\|u_{t}\right\|^2_{C_{L^2(\Omega)}}+\|h(t)\|^2.
\label{4.58-3}
\end{equation}

Then using $(\ref{1.2-3})-(\ref{1.5-3})$, we obtain
\begin{equation}
\begin{aligned}
\frac{d}{d t}\left(\|A^{\frac{\sigma}{2}}v_2(t)\|^2+|\varepsilon(t)|\|A^{\frac{1+\sigma}{2}}v_2(t)\|^2\right)& \leqslant C\left(\|A^{\frac{\sigma}{2}}v_2(t)\|^2+|\varepsilon(t)|\|A^{\frac{1+\sigma}{2}}v_2(t)\|^2\right)\\
&+C_g\left\|u_{t}\right\|^2_{C_{L^2(\Omega)}}+\|h(t)\|^2.
\end{aligned}
\label{4.59-3}
\end{equation}

Thanks to the Gronwall inequality, the summation formula for series and $h \in L_b^2\left(\mathbb{R} ; L^2(\Omega)\right)$, we arrive at
\begin{equation}
\begin{aligned}
\|A^{\frac{\sigma}{2}} v_2(t)\|^2+|\varepsilon(t)| \| A^{\frac{1+\sigma}{2}}&v_2(t) \|^2 \leqslant e^{C( t-\tau)}\left(1+C_g \int_\tau^t\left\|u_s\right\|_{C_{L^2(\Omega)}}^2 d s+\int_\tau^t\|h(s)\|^2 d s\right) \\
& \leqslant e^{C(t-\tau)}\left(1+C_g \int_\tau^t\left\|u_s\right\|_{C_{L^2(\Omega)}}^2 d s +(t-\tau)\|h(s)\|^2_{L_b^2(\mathbb{R};L^2(\Omega))}\right).
\end{aligned}
\label{4.60-3}
\end{equation}

Dividing $(\ref{4.22-3})$ by $e^{\beta t}$, then multiplying the resulting inequality by $C_ge^{C(t-\tau)}$, we deduce
\begin{align}
C_ge^{C(t-\tau)} \int_\tau^t\left\|u_s\right\|_{C_{L^2(\Omega)}}^2 d s & \leqslant C_ge^{C(t-\tau)}\int_\tau^t \frac{1}{1+\lambda_1 L}\left(\|\phi\|^2_{C_{L^2(\Omega)}}+\varepsilon_\tau\|\nabla \phi\|_{C_{L^2(\Omega)}}^2 \right)e^{-\beta_1(s-\tau)}ds \non\\
&+C_ge^{C(t-\tau)+\beta k} \int_\tau^t \frac{1}{\left.\delta\left(1+\lambda_1 L\right.\right)\left(1-e^{\left.-\beta_1\right.}\right.)}\|h(s)\|^2_{L_b^2(\mathbb{R};L^2(\Omega))} ds\non\\
&\leqslant e^{C(t-\tau)} \frac{C_g}{\beta_1\left(1+\lambda_1 L\right)} \left(\|\phi\|^2_{C_{L^2(\Omega)}}+\varepsilon_\tau\|\nabla \phi\|_{C_{L^2(\Omega)}}^2 \right)\non\\
&+e^{C(t-\tau)} \frac{C_g e^{\beta k}}{\delta\left(1+\lambda_1 L\right)\left(1-e^{\left.-\beta_1\right.})\right.}(t-\tau)\|h\|_{L_b^2\left(\mathbb R; L^2(\Omega)\right)}^2.
\label{4.61-3}
\end{align}

Inserting (\ref{4.61-3}) into (\ref{4.60-3}), we derive
\begin{equation}
\begin{aligned}
\|A^{\frac{\sigma}{2}} v_2(t)\|^2+|\varepsilon(t)| \| A^{\frac{1+\sigma}{2}}&v_2(t) \|^2 \leqslant e^{C( t-\tau)}\left(1+\frac{C_g}{\beta_1\left(1+\lambda_1 L\right)}\left(\|\phi\|^2_{C_{L^2(\Omega)}}+\varepsilon_\tau\|\nabla \phi\|_{C_{L^2(\Omega)}}^2 \right)\right)  \\
&\qquad\,\,\leqslant e^{C(t-\tau)}(t-\tau)\left(1+ \frac{C_g e^{\beta k}}{\delta\left(1+\lambda_1 L\right)\left(1-e^{\left.-\beta_1\right.})\right.}\right) \|h\|_{L_b^2\left(\mathbb R; L^2(\Omega)\right)}^2.
\end{aligned}
\label{4.62-3}
\end{equation}

Then putting $t+\theta$ instead of $t$ with $\theta \in [-k,0]$ in (\ref{4.62-3}) and by $(\ref{1.2-3})-(\ref{1.5-3})$, we obtain
\begin{align}
&\|A^{\frac{\sigma}{2}} v_{2,t}\|_{C_{L^2(\Omega)}}^2+|\varepsilon_t| \| A^{\frac{1+\sigma}{2}}v_{2,t} \|_{C_{L^2(\Omega)}}^2 \non\\
&\leqslant e^{C( t-k-\tau)}\left(1+\frac{C_g}{\beta_1\left(1+\lambda_1 L\right)} \left(\|\phi\|^2_{C_{L^2(\Omega)}}+\varepsilon_\tau
\|\nabla \phi\|_{C_{L^2(\Omega)}}^2 \right)\right) \non\\
&+e^{C(t-k-\tau)}(t-\tau)\left(1+ \frac{C_g e^{\beta k}}{\delta\left(1+\lambda_1 L\right)\left(1-e^{\left.-\beta_1\right.})\right.}\right) \|h\|_{L_b^2\left(\mathbb R; L^2(\Omega)\right)}^2.
\label{4.63-3}
\end{align}

Let $t=\tau+T_1$ in (\ref{4.63-3}), then by Lemmas \ref{lem2.11-3} and \ref{lem4.4-3}, we conclude there exists a constant $\widetilde{C}_1>0$ such that (\ref{4.52-3}) holds.  $\hfill$$\Box$

\begin{Lemma}\label{lem4.7-3}
Under the assumptions of Lemmas $\ref{lem4.3-3}$ and $\ref{lem4.4-3}$, assume $u(t)=\tilde v_1(t)+\tilde v_2(t)$ is a weak solution to problem $(\ref{1.1-3})$,
then there exists positive constants $C_{\xi}$ and $J_{\xi}$ such that the following inequalities hold
\begin{equation}
\int_s^t\left\|\nabla \tilde{v}_1(r)\right\|^2 d r \leqslant \xi(t-s)+C_{\xi},\, \text { for any } t \geqslant s \geqslant \tau,
\label{4.64-3}
\end{equation}
and
\begin{equation}
\|A^{\frac{1+\sigma}{2}} \tilde{v}_2(t)\|^2 \leqslant J_{\xi},\, \text { for any } t \geqslant \tau,
\label{4.65-3}
\end{equation}
where $\sigma$ is the same as in Lemma $\ref{lem4.4-3}$, $\xi>0$, $C_{\xi}$ and $J_{\xi}$ depend on $\xi$, $\|\phi\|^2_{C_{\mathcal H_t(\Omega)}}$ and $\|h\|_{L_b^2\left(\mathbb R; L^2(\Omega)\right)}^2$, and independent of $\tau$.
\end{Lemma}
$\mathbf{Proof.}$ From Lemma \ref{lem4.1-3}, we derive there exists a constant $\widetilde{C}_2>0$ depending on $\beta$, $\beta_1$, $L$, $\|\phi\|^2_{C_{\mathcal H_t(\Omega)}}$ and $\|h\|_{L_b^2\left(\mathbb R; L^2(\Omega)\right)}^2$ such that
\begin{equation}
\sup _{\tau \in \mathbb{R}} \sup _{\tau \leqslant t}\|U(t, \tau) u(\tau)\|_{C_{\mathcal H_t(\Omega)}}^2 \leqslant \widetilde{C}_2 .
\label{4.66-3}
\end{equation}

Then thanks to Lemma \ref{lem4.5-3}, (\ref{4.66-3}) and $u(t)=\tilde v_1(t)+\tilde v_2(t)$, we obtain there exists constants $\xi>0$ and $T_1>0$ such that
\begin{equation}
T_1 \geqslant \frac{1}{\widetilde{C}} \ln \frac{Z\left(\|\phi\|_{C_{\mathcal H_t(\Omega)}}^2\right)}{\xi},
\label{4.67-3}
\end{equation}
where $\widetilde{C}$ is the same as in (\ref{4.42-3}).

Fix a sufficiently large $T_1$ and let $\tilde v_1(t)=v_1(t)$ and $\tilde v_2(t)=v_2(t)$ in any interval $[\tau+(\eta-1)T_1, \tau+\eta T_1]$ with $\eta \in \mathbb{N^+}$, where $v_1(t)$ and $v_2(t)$ are weak solutions to problems (\ref{4.35-3}) and (\ref{4.36-3}), and their initial values are $v_1\left(\tau+(\eta-1) T_1\right)=U\left(\tau+(\eta-1) T_1, \tau\right) u(\tau)$ and $v_2\left(\tau+(\eta-1) T_1\right)=0$ respectively.

From (\ref{4.42-3}) and (\ref{4.66-3}), we conclude
\begin{equation}
\int_{\tau+(\eta-1) T_1}^{\left.\tau+ \eta T_1\right.}\left\|\tilde{v}_1(r)\right\|^2 d r \leqslant \widetilde{C}_2(B),
\label{4.68-3}
\end{equation}
where $B \subset C_{\mathcal H_{t}(\Omega)}$ is bounded, $\widetilde{C}_2$ depends on $B$ and independent of $\tau$, $\eta$ and $T_1$.

Hence, for any $\xi>0$, we can find a sufficiently large $T_1=T_1(\xi, B)$ and a constant $C_{\xi}>0$ such that (\ref{4.64-3}) holds.

In addition, if we fix the length of $T_1$, it follows that there exists a constant $J_{\xi}=J_{\xi}(T_1, B)$ yields (\ref{4.65-3}).  $\hfill$$\Box$

\begin{Corollary}\label{cor4.8-3}
Thanks to Lemmas $\ref{lem4.5-3}$ and $\ref{lem4.7-3}$, then there exists constants $Y_{\| \phi\|_{C_{\mathcal {H}_t(\Omega)}}^2}>0$ depending on $\|\phi\|_{C_{\mathcal {H}_t(\Omega)}}^2$ and $\widetilde C_3>0$ such that
\begin{equation}
\|\tilde{v}_1 (t)\|_{C_{\mathcal {H}_t(\Omega)}}^2 \leqslant Z\left(Y_{\| \phi\|_{C_{\mathcal {H}_t(\Omega)}}^2}\right):=\widetilde{C}_3,
\label{4.69-3}
\end{equation}
for any $\tau \leq t \in \mathbb{R}$.
\end{Corollary}

\begin{Lemma}\label{lem4.9-3}
Under the assumptions of Lemmas $\ref{lem4.3-3}$ and $\ref{lem4.4-3}$, for any bounded set $B \subset C_{\mathcal {H}_t(\Omega)}$, there exists a positive constant $W_{\|\phi\|_{C_{\mathcal {H}_t(\Omega)}}^2}>0$ depending on $\|\phi\|^2_{C_{\mathcal H_t(\Omega)}}$ such that
\begin{equation}
\|A^{\frac{\sigma}{2}} v_{2, t}\|_{C_{\mathcal {H}_t(\Omega)}}^2+|\varepsilon_t|\|A^{\frac{1+\sigma}{2}} v_{2, t}\|_{C_{\mathcal {H}_t(\Omega)}}^2 \leqslant W_{\|\phi\|^2_{C_{\mathcal {H}_t(\Omega)}}}
\label{4.70-3}
\end{equation}
for any $\tau \le t \in \mathbb R$ and $u(\tau)\in B$, where $\sigma$ is the same as in Lemma $\ref{lem4.4-3}$.
\end{Lemma}
$\mathbf{Proof.}$ Noticing $0<\sigma<\min \left\{\frac{1}{3}, \frac{n+2-(n-2) \gamma}{2}\right\}$ with $0<\gamma<\frac{n+2}{n-2}$ and $n \ge3$, then it is easily to deduce $\sigma < \frac{1+\sigma}{2}$.

By the Cauchy and Young inequalities, we obtain
\begin{equation}
2\left(h(t), A^{\sigma} v_2(t)\right) \leqslant \frac{1}{8}\|A^{\frac{1+\sigma}{2}} v_2(t)\|^2+8\|h(t)\|^2.
\label{4.71-3}
\end{equation}

From embedding $L^{\frac{2 n}{n-2}}(\Omega) \hookrightarrow L^{\frac{2 n \gamma}{n-2(1-\sigma)}}(\Omega)$ and $D(A^{\frac{1-\sigma}{2}}) \hookrightarrow L^{\frac{2 n}{n-2(1+\sigma)}}(\Omega)$, (\ref{1.10-3}) and the Young inequality, we conclude
\begin{align}
\left(f_1\left(v_1(t)\right), A^\sigma v_2(t)\right) & \leqslant C\int_{\Omega}\left(1+\left|v_1(t)\right|^{\gamma}\right)\left|A^{\sigma} v_2(t)\right| d x \non\\
& \leqslant C\left(1+\left\|v_1(t)\right\|^{\gamma}_{L^{\frac{2 n \gamma}{n+2-2\sigma}}(\Omega)}\right) \|A^{\sigma} v_2(t)\|_{L^{\frac{2 n}{n-2+2\sigma}}(\Omega)} \non\\
& \leqslant C\left(1+\left\|v_1(t)\right\|^{\gamma}_{L^{\frac{2n}{n-2}}(\Omega)}\right) \left\|A^{\sigma} v_2(t)\right\|_{L^{\frac{2 n}{n-2+2\sigma}}(\Omega)} \non\\
& \leqslant C\left(1+\left\|v_1(t)\right\|^{\gamma}_{L^{\frac{2n}{n-2}}(\Omega)}\right) \|A^{\frac{1+\sigma}{2}} v_2(t)\| \non\\
& \leqslant C\left(1+\left\|v_1(t)\right\|_{L^{\frac{2 n}{n-2}}(\Omega)}^{2 \gamma}\right)+\frac{1}{16}\|A^{\frac{1+\sigma}{2}} v_2(t)\|^2.
\label{4.72-3}
\end{align}

Using (\ref{1.7-3}), we derive
\begin{equation}
\begin{aligned}
|(f(u(t))-f\left(v_1(t)), A^{\sigma} v_2(t))|\right. & \leqslant C \int_{\Omega}\left|(f^{\prime}\left((1-\mu)u(t))+\mu v_1(t)\right)\right||u(t)-v_1(t)|\,| A^{\sigma} v_2(t)|d x \\
& \leqslant C \int_{\Omega}\left(1+|u(t)|^{\frac{4}{n-2}}+\left|v_1(t)\right|^{\frac{4}{n-2}}\right)\left|v_2(t)\right|\left|A^{\sigma} v_2(t)\right| d x,
\end{aligned}
\label{4.73-3}
\end{equation}
where $0<\mu<1$.

Noticing $u(t)=v_1(t)+v_2(t)$, it follows that
\begin{equation}
\int_{\Omega}|u(t)|^{\frac{4}{n-2}}\left|v_2(t)\right||A^{\sigma} v_2(t)| d x \leq C \int_{\Omega}\left(\left|v_1(t)\right|^{\frac{4}{n-2}}+\left|v_2(t)\right|^{\frac{4}{n-2}}\right)\left|v_2(t)\right|\left|A^{\sigma} v_2(t)\right| d x.
\label{4.74-3}
\end{equation}

Hence, from Lemma \ref{lem2.11-3}, the Cauchy and Young inequalities, we obtain there exists positive constants $\widetilde C_4$ and $\widetilde C_5$ such that
\begin{equation}
\int_{\Omega}\left|v_2(t)\right|\left|A^{\sigma} v_2(t)\right| d x \leqslant \widetilde{C}_4+\widetilde{C}_5\|A^{\frac{1+{\sigma}}{2}} v_2(t)\|^2.
\label{4.75-3}
\end{equation}

Conducting similar calculations to (\ref{4.72-3}), then by Corollary \ref{cor4.8-3}, we deduce
\begin{align}
\int_{\Omega}|v_1(t)|^{\frac{4}{n-2}}|v_2(t)||A^{\sigma} v_2(t)| d x & \leqslant C\left\|v_1(t)\right\|_{L^{\frac{2n}{n-2}}(\Omega)}^{\frac{4}{n-2}}\left\|v_2(t)\right\|_{L^{\frac{2 n}{n-2(1+\sigma)}}(\Omega)}\left\|A^{\sigma} v_2(t)\right\|_{L^{\frac{2 n}{n-2(1-\sigma)}}(\Omega)}\non\\
&\leqslant C\|A^{\frac{1}{2}} v_1(t)\|^{\frac{4}{n-2}}\|A^{\frac{1+\sigma}{2}} v_2(t)\|^2 \non\\
& \leqslant \frac{1}{16}\|A^{\frac{1+\sigma}{2}} v_2(t)\|^2+C \widetilde{C}_3\left\|\nabla v_1(t)\right\|^2\|A^{\frac{1+\sigma}{2}} v_2(t)\|^2.
\label{4.76-3}
\end{align}

Thanks to Lemma $\ref{lem4.4-3}$ and (\ref{4.65-3}), by similar calculations as those in Lemma 4.7 in \cite{sy.3}, we derive
\begin{equation}
\int_{\Omega}\left|v_2(t)\right|^{\frac{4}{n-2}}\left|v_2(t)\right||A^{\sigma} v_2(t)| d x \leq CJ_\xi^{\frac{4}{n-2}}+\frac{1}{16}\|A^{\frac{1+\sigma}{2}} v_2(t)\|^2
\label{4.77-3}
\end{equation}
and
\begin{equation}
\int_{\Omega}\left|v_1(t)\right|^{\frac{4}{n-2}}\left|v_2(t)\right||A^{\sigma} v_2(t)| d x \leq C\|\nabla v_1(t)\|^{\frac{4}{n-2}}\|A^{\frac{1+\sigma}{2}} v_2(t)\|^2.
\label{4.78-3}
\end{equation}

Inserting $(\ref{4.74-3})-(\ref{4.78-3})$ into (\ref{4.73-3}), we arrive at
\begin{equation}
\begin{aligned}
|(f(u(t))-f(v_1(t)), A^{\sigma} v_2(t))|& \leqslant (\frac{1}{8}+\widetilde C_5)\|A^{\frac{1+\sigma}{2}} v_2(t)\|^2+C\|\nabla v_1(t)\|^{\frac{4}{n-2}}\|A^{\frac{1+\sigma}{2}} v_2(t)\|^2\\
& +CJ_\xi^{\frac{4}{n-2}}+C \widetilde{C}_3\left\|\nabla v_1(t)\right\|^2\|A^{\frac{1+\sigma}{2}} v_2(t)\|^2+\widetilde{C}_4.
\end{aligned}
\label{4.79-3}
\end{equation}

From $(\ref{1.11-3})-(\ref{1.13-3})$ and the Young inequality, we deduce
\begin{equation}
2(g(t, u_{t}), A^{\sigma}v_2(t)) \leq 2 C_{g} \|u_{t}\|^{2}_{C_{L^{2}(\Omega)}}+\|A^{\sigma}v_2(t)\|^{2}.
\label{4.80-3}
\end{equation}

Inserting (\ref{4.71-3}), (\ref{4.72-3}), (\ref{4.79-3}) and (\ref{4.80-3}) into (\ref{4.53-3}), we conclude there exists a constant $\widetilde C_6>0$ such that
\begin{equation}
\begin{aligned}
& \frac{d}{d t}\left(\|A^{\frac{\sigma}{2}}v_2(t)\|^2+\varepsilon(t)\|A^{\frac{1+\sigma}{2}}v_2(t)\|^2\right)+\left(2 a(l(u))-\varepsilon^{\prime}(t)\right)\|A^{\frac{1+\sigma}{2}} v_2(t)\|^2 \\
& \leqslant CJ_{\xi}^{\frac{4}{n-2}}+C\left(1+\left\|v_1(t)\right\|_{L^{\frac{2n}{n-2}}(\Omega)}^{2\gamma}\right)+C \widetilde C_6 \|\nabla v_1(t)\|^2\|A^{\frac{1+\sigma}{2}}v_2(t)\|^2+8\|h(t)\|^2.
\end{aligned}
\label{4.81-3}
\end{equation}

Then using assumptions $(\ref{1.2-3})-(\ref{1.5-3})$, $(\ref{4.42-3})$ and Lemma \ref{lem4.7-3}, we deduce there exists constants satisfying $C-\widetilde{C}_7\left\|\nabla v_1(t)\right\|^2>0$, $\widetilde{C}_7=C\widetilde{C}_6$ and $\widetilde C_8>C(1+J_{\xi}^{\frac{4}{n-2}})$ such that
\begin{equation}
\begin{aligned}
&\frac{d}{d t}\left(\|A^{\frac{\sigma}{2}}v_2(t)\|^2+|\varepsilon(t)|\|A^{\frac{1+\sigma}{2}}v_2(t)\|^2\right)+(C-\widetilde{C}_7\|\nabla v_1(t) \|^2)(\|A^{\frac{\sigma}{2}} v_2(t)\|^2\\
&+ |\varepsilon(t)|\|A^{\frac{1+\sigma}{2}} v_2(t)\|^2) \leqslant \widetilde{C}_8+C_g\left\|u_t\right\|_{C_{L^2(\Omega)}}^2+8\|h(t)\|^2.
\end{aligned}
\label{4.82-3}
\end{equation}

By the Gronwall inequality, we conclude
\begin{align}
& \|A^{\frac{\sigma}{2}} v_2(t)\|^2+|\varepsilon(t)|\|A^{\frac{1+\sigma}{2}}v_2(t)\|^2\non\\
& \leqslant e^{-\int_{\tau+T_2}^t(C-\widetilde{C}_7\|\nabla v_1(s) \|^2)d s}\left(\|A^{\frac{\sigma}{2}} v_2(\tau+T_2)\|^2+\varepsilon(\tau+T_2)\|A^{\frac{1+\sigma}{2}}v_2(\tau+T_2)\|^2\right) \non\\
& +\widetilde{C}_8 \int_{\tau+T_2}^t e^{\int_t^s(C-\widetilde{C}_7\|\nabla v_1(y)\|^2) d y} d s+C_g \int_{\tau+T_2}^t e^{\int_t^s(C-\widetilde{C}_7\|\nabla v_1(y) \|^2) d y}\|u_s\|_{C_{L^2(\Omega)}}^2 d s \non\\
& +8 \int_{t+T_2}^t\|h(s)\|^2 e^{\int_t^s(C-\widetilde{C}_7\|\nabla v_1(y) \|^2) d y} d s.
\label{4.83-3}
\end{align}

Suppose $\xi<\frac{C}{2\widetilde{C}_7}$, then by $(\ref{4.22-3})$ and $(\ref{4.64-3})$, we derive
\begin{align}
\int_{\tau+T_2}^{t}e^{\int_t^s(C-\widetilde{C}_7\|\nabla v_1(y) \|^2)d y} \|u_s\|_{C_{L^2(\Omega)}}^2 d s & \leqslant e^{\widetilde {C}_7 C_{\xi}}\int_{\tau+T_2}^{t}\|u_s\|_{C_{L^2(\Omega)}}^2 e^{-\frac{c}{2}(t-s)}ds\non\\
&\leqslant \frac{2}{C(1+\lambda_1 L)} e^{\widetilde {C}_7 C_{\xi}}\left(\|\phi\|^2_{C_{L^2(\Omega)}}+\varepsilon_\tau\|\nabla \phi\|_{C_{L^2(\Omega)}}^2 \right)\non\\
& +\frac{1}{C \delta\left(1+\lambda_1 L\right)\left(1-e^{\left.-\beta_1\right)}\right.} e^{\widetilde {C}_7 C_{\xi}+\beta k}\|h\|_{L_b^2(\mathbb{R} ; L^2(\Omega))}^2.
\label{4.84-3}
\end{align}

Through similar calculations as those in (\ref{4.84-3}), we obtain
\begin{equation}
\begin{aligned}
\int_{\tau+T_2}^{t}\|h(s)\|^2e^{\int_t^s(C-\widetilde{C}_3\|\nabla v_1(y) \|^2)d y} d s & \leqslant e^{\widetilde {C}_7 C_{\xi}} e^{-\frac{c}{2}t}\int_{\tau+T_2}^{t} e^{\frac{c}{2}s}\|h(s)\|^2ds\\
&\leqslant \frac{e^{\widetilde {C}_7 C_{\xi}}}{1-e^{-\frac{c}{2}}}\|h\|_{L_b^2(\mathbb{R} ; L^2(\Omega))}^2.\\
\end{aligned}
\label{4.85-3}
\end{equation}

Inserting (\ref{4.84-3}) and (\ref{4.85-3}) into (\ref{4.83-3}), we arrive at
\begin{equation}
\begin{aligned}
 \|A^{\frac{\sigma}{2}} v_2(t)\|^2+|\varepsilon(t)|\|A^{\frac{1+\sigma}{2}}v_2(t)\|^2& \leqslant C\left(\|A^{\frac{\sigma}{2}} v_2(\tau+T_2)\|^2+\varepsilon(\tau+T_2)\|A^{\frac{1+\sigma}{2}}v_2(\tau+T_2)\|^2\right) \\
& +C\left(\|\phi\|^2_{C_{L^2(\Omega)}}+\varepsilon_\tau\|\nabla \phi\|_{C_{L^2(\Omega)}}^2 \right)+C\|h\|_{L_b^2(\mathbb{R} ; L^2(\Omega))}^2.
\end{aligned}
\label{4.86-3}
\end{equation}

Putting $t+\theta$ instead of $t$ with $\theta \in [-k,0]$ in (\ref{4.86-3}) yields
\begin{equation}
\begin{aligned}
 \|A^{\frac{\sigma}{2}} v_{2,t}\|^2+|\varepsilon_t|\|A^{\frac{1+\sigma}{2}}v_{2,t}\|^2& \leqslant C\left(\|A^{\frac{\sigma}{2}} v_2(\tau+T_2)\|^2+\varepsilon(\tau+T_2)\|A^{\frac{1+\sigma}{2}}v_2(\tau+T_2)\|^2\right) \\
& +C\left(\|\phi\|^2_{C_{L^2(\Omega)}}+\varepsilon_\tau\|\nabla \phi\|_{C_{L^2(\Omega)}}^2 \right)+C\|h\|_{L_b^2(\mathbb{R} ; L^2(\Omega))}^2.
\end{aligned}
\label{4.87-3}
\end{equation}

Therefore, if $B \subset C_{\mathcal H_t(\Omega)}$ is bounded and $\phi \subset B$, then there exists a constant $W_{\|\phi\|^2_{C_{\mathcal {H}_t(\Omega)}}}>0$ depending on $\|\phi\|^2_{C_{\mathcal H_t(\Omega)}}$ such that (\ref{4.70-3}) holds. $\hfill$$\Box$

\begin{Lemma}\label{lem4.10-3}
Under the assumptions of Lemmas $\ref{lem4.3-3}$ and $\ref{lem4.4-3}$, if $v_2(t)$ is a weak solution to problem $(\ref{4.36-3})$, then it satisfies $\partial_t v_2 \in L^2(\tau, t ;\mathcal H_t(\Omega))$ for any $\tau \le t \in \mathbb R$.
\end{Lemma}
$\mathbf{Proof.}$ Taking $L^2(\Omega)$-inner product between $\partial_t v_1$ and $(\ref{4.35-3})_1$, by $(\ref{1.2-3})-(\ref{1.5-3})$, (\ref{1.9-3}), the Cauchy and Young inequalities, we obtain
\begin{equation}
\left\|\partial_t v_1\right\|^2+2 \varepsilon(t)\left\|\partial_t \nabla v_1\right\|^2+m \frac{d}{d t}\left\|\nabla v_1\right\|^2 \leqslant\|f_0 v_1 \|^2.
\label{4.88-3}
\end{equation}

Then integrating (\ref{4.88-3}) from $\tau$ to $t$ and using $(\ref{1.2-3})-(\ref{1.5-3})$, it follows that $\partial_t v_1 \in L^2(\tau, t ;\mathcal H_t(\Omega))$.

Furthermore, since $\partial_t v_2=\partial_t u-\partial_t v_1$ and noticing $\partial_t u \in L^2(\tau, t ;\mathcal H_t(\Omega))$, then $\partial_t v_2 \in L^2(\tau, t ;\mathcal H_t(\Omega))$ follows directly. $\hfill$$\Box$

\begin{Lemma}\label{lem4.11-3} Under the assumptions of Lemmas $\ref{lem4.3-3}$ and $\ref{lem4.4-3}$, then the process $\{U(t, \tau)\}_{t \geq \tau}$ of weak solutions to problem $(\ref{1.1-3})$ is pullback $\mathcal D_{C_{\mathcal H_t(\Omega)}}$-$\omega$-limit compact in the time-dependent space $C_{\mathcal H_t(\Omega)}$.
\end{Lemma}
$\mathbf{Proof.}$ Assume the sequences $\left\{u^k(t)\right\}_{k\in \mathbb N^{+}}$, $\left\{v_1^k(t)\right\}_{k \in \mathbb N^{+}}$ and $\left\{v_2^k(t)\right\}_{k \in \mathbb N^{+}}$ are weak solutions to problems (\ref{1.1-3}), (\ref{4.35-3}) and (\ref{4.36-3}) respectively, which satisfy $u^k(t)=v_1^k(t)+v_2^k(t)$.

From Lemma \ref{lem4.9-3}, we deduce $\left\{v_{2,t}^k\right\}_{k \in \mathbb N^{+}}=\left\{v_2^k(t+\theta)\right\}_{k \in \mathbb N^{+}}$ is bounded in $C_{\mathcal{H}_t(\Omega), \sigma}$, then we conclude $v_{2}^k \in L^{\infty}(t-k, t; H^{1+\sigma}(\Omega))$, where $H^{1+\sigma}(\Omega)$ is equipped with the norm $\|A^\frac{\sigma}{2}\cdot\|+\|A^\frac{1+\sigma}{2}\cdot\|$ with $A=-\Delta$.

Integrating (\ref{4.70-3}) from $t-k$ to $t$, we arrive at
\begin{equation}
\int_{t-k}^t\left\|v_{2, t}\right\|_{\mathcal H_t(\Omega)}^2 d t \leq \int_{t-k}^t \|v_{2, t}\|^2_{C_{\mathcal H_t(\Omega)}}dt \leq C  J_{\|\phi\|^2_{C_{\mathcal {H}_t(\Omega)}}},
\label{4.89-3}
\end{equation}
where the constant $J_{\|\phi\|^2_{C_{\mathcal {H}_t(\Omega)}}}>0$ depends on $\|\phi\|^2_{C_{\mathcal H_t(\Omega)}} \subset B$.

Thanks to Lemma \ref{lem4.10-3}, we obtain $\partial_t v_2^k \in L^2(t-k, t ;\mathcal H_t(\Omega))$. Then we deduce $v_2^k \in L^{\infty}(t-k, t; H^{1+\sigma}(\Omega)) \cap W^{1, 2}(t-k, t ;\mathcal H_t(\Omega))$. Hence, we deduce $\left\{v_{2,t}^k\right\}_{k \in \mathbb N^{+}}$ is relatively compact in $C([t-k,t];\mathcal H_t(\Omega))$.

Consequently, by Lemma \ref{lem2.9-3}, the pullback $\mathcal D_{C_{\mathcal H_t(\Omega)}}$-$\omega$-limit compactness of the process holds directly. $\hfill$$\Box$

\begin{Theorem}\label{th4.12-3}
Under the assumptions of Lemmas $\ref{lem4.3-3}$ and $\ref{lem4.4-3}$, then there exists a nonempty set $\widehat{D}_{0}$ satisfying Definition $\ref{def4.2-3}$ such that the process $\{U(t, \tau)\}_{t \geq \tau}$ of weak solutions to problem $(\ref{1.1-3})$ processes a unique pullback $\mathcal D_{C_{\mathcal H_t(\Omega)}}$-attractor $\mathcal A$ in the time-dependent space $C_{\mathcal H_t(\Omega)}$.
\end{Theorem}
$\mathbf{Proof.}$ From Lemmas $\ref{lem2.15-3}$, $\ref{lem4.1-3}$, $\ref{lem4.3-3}$, $\ref{lem4.5-3}-\ref{lem4.7-3}$, $\ref{lem4.9-3}$ and $\ref{lem4.11-3}$, we deduce there exists a nonempty set $\widehat{D}_{0}$ satisfying Definition $\ref{def4.2-3}$, then by Definition $\ref{def2.14-3}$, it follows there exists a unique pullback $\mathcal D_{C_{\mathcal H_t(\Omega)}}$-attractor $\mathcal A$ in $C_{\mathcal H_t(\Omega)}$, which attracts every bounded subset of $C_{\mathcal H_t(\Omega)}$. $\hfill$$\Box$

\section{\large Regularity of pullback attractors}
In Lemmas \ref{lem4.6-3}, \ref{lem4.7-3} and \ref{lem4.9-3}, we perform some prior estimates for the regularity of weak solutions to problem $(\ref{1.1-3})$ in the time-dependent space $C_{\mathcal{H}_t^1(\Omega),\sigma}$. In this section, we will further discuss the regularity in the time-dependent space $C_{\mathcal{H}_t^1(\Omega)}$ by decomposing the solution $u$ into two parts and conducting a series of calculations and estimates.
\begin{Theorem}\label{th5.1-3}
Under the assumptions of Lemmas $\ref{lem4.3-3}$ and $\ref{lem4.4-3}$, the pullback $\mathcal D_{C_{\mathcal H_t(\Omega)}}$-attractor $\mathcal A$ is bounded in $C_{\mathcal{H}_t^1(\Omega)}$.
\end{Theorem}
$\mathbf{Proof.}$ Noticing that the regularity of $C_{\mathcal{H}_t^1(\Omega),\sigma}$ is higher than $C_{\mathcal{H}_t(\Omega)}$, and from Lemmas \ref{lem4.6-3}, \ref{lem4.7-3} and \ref{lem4.9-3}, we derive there exist $\tilde{h} \in L_{b}^{2}(\mathbb R ; C_{\mathcal{H}_t(\Omega)})$ and constant $r_1>0$ such that
\begin{equation}
\|h-\tilde{h}\|<r_1^2.
\label{5.1-3}
\end{equation}

Then assume $\phi \in \mathcal A$, $u^1(t)=U_1(t,\tau+\theta ){u(\tau+\theta)}$ and  $u^2(t)=U_2(t,\tau+\theta ){u(\tau+\theta)}$ with $u(t)=U(t, \tau+\theta) u(\tau+\theta)=U_{1}(t, \tau+\theta) u(\tau+\theta)+U_{2}(t, \tau+\theta) u(\tau+\theta)$ being a weak solution to  problem $(\ref{1.1-3})$.

Furthermore, assume $u^1$ satisfies
\begin{equation}
\left\{\begin{array}{ll}
\partial_{t} u^1-\varepsilon(t)\partial_{t} \Delta u^1-a(l(u)) \Delta u^1=h-\tilde{h} & \text { in } \Omega \times(\tau, \infty), \\
u^1(x,t)=0 & \text { on } \partial \Omega \times(\tau, \infty), \\
u^1(x, \tau+\theta)=\phi(x,\theta), &\,\, x \in \Omega,\, \theta \in[-k, 0],
\label{5.2-3}
\end{array}\right.
\end{equation}
and $u^2$ satisfies
\begin{equation}
\left\{\begin{array}{ll}
\partial_{t} u^2-\varepsilon(t)\partial_{t} \Delta u^2-a(l(u)) \Delta u^2=f(u)+g(t,u_{t})+\tilde{h} & \text { in } \Omega \times(\tau, \infty), \\
u^2(x,t)=0 & \text { on } \partial \Omega \times(\tau, \infty), \\
u^2(x, \tau+\theta)=0, &\,\, x \in \Omega,\, \theta \in[-k, 0].
\label{5.3-3}
\end{array}\right.
\end{equation}

Taking $L^2(\Omega)$-inner product between $- \Delta u^1$ and $(\ref{5.2-3})_{1}$ and using (\ref{5.1-3}), we obtain
\begin{equation}
\frac{d}{d t}(\|\nabla u^1\|^2+\varepsilon(t)\|\Delta u^1\|^2)+(2 a(l(u))-\varepsilon^{\prime}(t)-1)\|\Delta u^1\|^2\leq r_1^2.
\label{5.4-3}
\end{equation}

Suppose $\tilde{r}_1>0$, then by the Poincar\'{e} inequality, we deduce
\begin{equation}
\tilde r_1(\|\nabla u^1\|^2+|\varepsilon(t)|\|\Delta u^1\|^2)\leq \tilde r_1(\lambda_1^{-1}+|\varepsilon(t))|\|\Delta u^1\|^2.
\label{5.5-3}
\end{equation}

From $(\ref{1.2-3})-(\ref{1.5-3})$ and noticing $m>\frac{3}{2}+\frac{L}{2}+\frac{1}{4 \lambda_1}$, we derive when $\varepsilon(t)$ is decreasing, $$\min\left\{\frac{2a(l(u))-\varepsilon^{\prime}(t)-1}{\lambda_1^{-1}+|\varepsilon(t)|} \right\}=\frac{2+L+(2\lambda_1)^{-1}}{\lambda_1^{-1}+L} $$
and when $\varepsilon(t)$ is increasing, $$\min\left\{\frac{2a(l(u))-\varepsilon^{\prime}(t)-1}{\lambda_1^{-1}+|\varepsilon(t)|} \right\}=\frac{2+2L+(2\lambda_1)^{-1}}{\lambda_1^{-1}+L}.$$
Hence, if $\tilde r_1>0$ further satisfies
$0<\tilde{r}_1 \leqslant \frac{2+L+(2\lambda_1)^{-1}}{\lambda_1^{-1}+L}$, it follows
\begin{equation}
\frac{d}{d t}(\|\nabla u^1\|^2+|\varepsilon(t)|\|\Delta u^1\|^2)+\tilde r_1(\|\nabla u^1\|^2+|\varepsilon(t)|\|\Delta u^1\|^2)\le r_1^2,
\label{5.6-3}
\end{equation}
where $\alpha$ is the same as in $(\ref{1.2-3})$.

Let $B_{1}(t)=\|\nabla u^1\|^2+|\varepsilon(t)|\|\Delta u^1\|^2$, we derive
\begin{equation}
\frac{d}{d t} B_{1}(t)+ \tilde r_1 B_{1}(t) \leq r_1^2.
\label{5.7-3}
\end{equation}

By the Gronwall inequality and putting $t+\theta$ instead of $t$ in the resulting equation, it follows
\begin{equation}
\|u^1\|_{C_{\mathcal{H}_{t}^{1}(\Omega)}}^{2} \leq K_1,
\label{5.8-3}
\end{equation}
where $K_1=e^{-\tilde r_1(t+k-\tau)}\left\|\phi\right\|_{C_{\mathcal{H}_{t}^{1}(\Omega)}}^{2}+\frac{r_1^2}{\tilde r_1}$ is a positive constant.

Then taking $L^2(\Omega)$-inner product between $- \Delta u^2$ and $(\ref{5.3-3})_{1}$, we conclude
\begin{equation}
\frac{d}{d t}\left(\|\nabla u^2\|^{2}+\varepsilon(t)\|\Delta u^2\|^{2}\right)+\left(2 a(l(u))-\varepsilon^{\prime}(t)\right)\|\Delta u^2\|^{2}=2(f(u)+g(t,u_{t})+\tilde h,-\Delta u^2).
\label{5.9-3}
\end{equation}

Using the Young inequality, assumptions $(\ref{1.2-3})-(\ref{1.6-3})$ and $(\ref{1.11-3})-(\ref{1.13-3})$, we arrive at
\begin{equation}
2(f(u),-\Delta u^2) \leq 2 \lambda_1^2 \|u\|^{2}+\frac{1}{2}\|\Delta u^2\|^{2},
\label{5.10-3}
\end{equation}
\begin{equation}
2(g(t, u_{t}),-\Delta u^2) \leq 4C_{g}\|u_{t}\|_{C_{L^{2}(\Omega)}}^{2}+\frac{1}{4}\|\Delta u^2\|^{2}
\label{5.11-3}
\end{equation}
and
\begin{equation}
2(\tilde h,-\Delta u^2) \leq 4\|\tilde h\|^{2}+\frac{1}{4}\|\Delta u^2\|^{2}.
\label{5.12-3}
\end{equation}

Inserting (\ref{5.10-3})$-$(\ref{5.12-3}) into $(\ref{5.9-3})$ and by the Poincar\'{e} inequality, we obtain
\begin{equation}
\begin{aligned}
& \frac{d}{d t}(\|\nabla u^2\|^2+\varepsilon(t)\|\Delta u^2\|^2)+\left(2 a (l(u))-\varepsilon^{\prime}(t)-1\right)\|\Delta u^2\|^2 \\
& \leqslant 2 \lambda_1^2\|u\|^2+4C_{g}\|u_{t}\|_{C_{L^{2}(\Omega)}}^{2}+4\|\tilde{h}\|^2.
\end{aligned}
\label{5.13-3}
\end{equation}

Similarly, from Lemma \ref{lem4.3-3}, the Poincar\'{e} inequality and $(\ref{1.2-3})-(\ref{1.5-3})$, we deduce there exists a constant $0<\tilde{r}_2 \leqslant \frac{2+L+(2\lambda_1)^{-1}}{\lambda_1^{-1}+L}$ such that
\begin{equation}
\frac{d}{d t}(\|\nabla u^2\|^2+|\varepsilon(t)|\|\Delta u^2\|^2)+\tilde r_2(\|\nabla u^2\|^2+|\varepsilon(t)|\|\Delta u^2\|^2)\le C_{\lambda_1, R(t), C_g,\tilde h}R^2(t),
\label{5.14-3}
\end{equation}
where $R(t)$ is the same as in Lemma $\ref{lem4.3-3}$.

Then taking $B_{2}(t)=\|\nabla u^2\|^2+|\varepsilon(t)|\|\Delta u^2\|^2$, we deduce
\begin{equation}
\frac{d}{d t} B_{2}(t)+ \tilde r_2 B_{2}(t) \leq C_{\lambda_1, R(t), C_g,\tilde h}R^2(t).
\label{5.15-3}
\end{equation}

Therefore, it follows
\begin{equation}
\|u^2\|_{C_{\mathcal{H}_{t}^{1}(\Omega)}}^{2} \leq K_2,
\label{5.16-3}
\end{equation}
where $K_2=e^{-\tilde r_2(t+k-\tau)}\left\|\phi\right\|_{C_{\mathcal{H}_{t}^{1}(\Omega)}}^{2}+C_{\lambda_1, R(t), C_g,\tilde h}e^{-\tilde r_2(t-k)}\int_{\tau}^{t}e^{\tilde r_2s}R^2(s)ds$.

Combining $(\ref{5.8-3})$ with $(\ref{5.16-3})$, we conclude
\begin{equation}
\|u\|_{C_{\mathcal{H}_{t}^{1}(\Omega)}}^{2} \leq K_1+K_2=\bar K,
\label{5.17-3}
\end{equation}
which gives
\begin{equation}
\lim _{\tau \rightarrow-\infty} \operatorname{dist}(\mathcal A, \bar{\mathscr B}_{1, C_{\mathcal{H}_{t}^{1}(\Omega)}}(\bar K))=0,
\label{5.18-3}
\end{equation}
where
$$
\bar{\mathscr B}_{1, C_{\mathcal{H}_t^1(\Omega)}}\left(\bar K\right)=\left\{u(t) \in \bar{\mathscr B}_{1, C_{\mathcal{H}_t^1(\Omega)}}:\|u(t)\|_{C_{\mathcal{H}_t^1(\Omega)}}^2 \leq \bar K\right\}.
$$

Hence, we derive that $\mathcal{A} \subseteq \bar{\mathscr B}_{1, C_{\mathcal{H}_t^1(\Omega)}}\left(\bar K\right)$, which completes the proof. $\hfill$$\Box$

$\mathbf{Acknowledgment}$

This paper was supported by the China Scholarship Council with number 202206630048, the National Natural Science Foundation of China with contract number 12171082, the fundamental research funds for the central universities with contract numbers $2232022G$-$13$, $2232023G$-$13$ and a grant from science and technology commission of Shanghai municipality.

$\mathbf{Conflict\,\,of\,\,interest\,\,statement}$

The authors have no conflict of interest.

\end{document}